\documentclass[a4paper, 11pt, final]{amsart}

\makeatletter

\numberwithin{equation}{section}
\numberwithin{figure}{section}
\usepackage[a4paper]{geometry}
\usepackage{amsmath,amsfonts,amssymb,amsthm,epsfig}
\usepackage{amsmath}
\usepackage{color,xcolor}
\usepackage[final,linkcolor = blue,citecolor = blue,colorlinks=true]{hyperref}
\usepackage[leqno]{amsmath}

\usepackage{amscd}
\usepackage{mathrsfs}
\usepackage{}
\usepackage[all]{xy}
\usepackage[OT1]{fontenc}
\usepackage[latin1]{inputenc}
\usepackage{amssymb,latexsym}
\usepackage{type1cm,ae}
\usepackage[notref,notcite]{showkeys}
\usepackage{graphicx}
\usepackage{xspace}
\usepackage{setspace}
\usepackage{enumerate}
\usepackage[english]{babel}
\usepackage{bbm}
\usepackage{esint}
\usepackage{soul}

\soulregister{\cite}7
\soulregister{\citep}7
\soulregister{\citet}7
\soulregister{\ref}7
\soulregister{\pageref}7
\sethlcolor{yellow}

\theoremstyle{definition}
\newtheorem{theorem}{Theorem}[section]

\newtheorem{lemma}[theorem]{Lemma}

\newtheorem{definition}[theorem]{Definition}
\newtheorem{remark}[theorem]{Remark}

\numberwithin{equation}{section}

\newcommand*{\Wert}{\mathord{\mbox{|\kern-1.5pt|\kern-1.5pt|}}}

\renewcommand{\P}{\mathbb{P}}	

\DeclareMathOperator{\supp}{supp}

\def\rr{{\Bbb R}}
\def\rz{{{\rr}^n}}

\def\nn{{\Bbb N}}
\def\cc{{\Bbb C}}

\def\L{\mathcal{L}}

\def\Ez{\mathcal{E}}

\def\P{\mathcal{P}}

\def\H{\mathcal{H}}

\def\rdm{\rr^{n+1}}

\def\fz{\infty}

\def\supp{{\rm{\ supp\ }}}

\def\ez{\epsilon}

\def\supp{{\rm supp}}

\def\l{\left}
\def\r{\right}

\def\XXint#1#2#3{{\setbox0=\hbox{$#1{#2#3}{\int}$}
		\vcenter{\hbox{$#2#3$}}\kern-.5\wd0}}

\makeatother

\usepackage{babel}

\begin{document}
\raggedbottom
\allowdisplaybreaks

\title[Gaussian estimates for higher-order parabolic equations]{Gaussian estimates for fundamental solutions of higher-order parabolic equations with time-independent coefficients}
	
	\author[Guoming Zhang]{Guoming Zhang}
	
	\address[Guoming Zhang]{College of Mathematics and System Science, Shandong University of Science and Technology, Qingdao, 266590, Shandong, People's Republic of China}
	\email{ zhangguoming256@163.com}

	

	\thanks{$^*$Corresponding author: Guoming Zhang}
	\thanks{The author is supported by the National Natural Science Foundation of Shandong Province (No. ~ZR2023QA124).}

	\thanks{}
	
	\date{}
	
	\begin{abstract} We study the De Giorgi-Moser-Nash estimates of higher-order parabolic equations in divergence form with complex-valued, measurable, bounded, uniformly elliptic (in the sense of G$\mathring{a}$rding inequality) and time-independent coefficients. We also obtain
Gaussian upper bounds and H\"{o}lder regularity estimates for the fundamental solutions of this class of parabolic equations.

	\end{abstract}

	\maketitle
	\tableofcontents
	
\section{Introduction}

In \cite{DA, FS}, Gaussian upper and lower bounds for the fundamental solutions of second order parabolic equations in divergence form with real and bounded coefficients were established. Later, under the assumption that the complex-valued and time-independent coefficients are a small $L^{\infty}$ perturbation of real coefficients, Auscher \cite{A} showed that Gaussian upper bound holds for the fundamental solutions of second order parabolic equations of this type. Subsequently, Hofmann and Kim \cite{HK} proved that Auscher's result remain valid for a parabolic system of second order with time-dependent coefficients if the system is a small perturbation of a diagonal system. In particular, inspired by the result of Auscher, McIntosh and Tchamitchian \cite{AMT}, Kim \cite{K} showed that the fundamental solutions of second order parabolic systems with time-independent and complex-valued coefficients in $\rr^{2}$ have Gaussian upper bounds in the absence of the small perturbation assumption. More importantly, Kim constructed the H\"{o}lder continuity estimates for weak solutions of a second-order parabolic system with time-dependent coefficients if weak solutions of the corresponding elliptic system satisfy the H\"{o}lder estimates at every scale. 

In 2000, Auscher and Qafsaoui \cite{AQ} extended the results in \cite{A} for second order elliptic operators to elliptic operators of any order, in which they found a criterion to decide on whether upper bounds and H\"{o}lder regularity hold for the heat kernel of divergence form elliptic operators (systems) of higher order with complex coeffcients on $\rz.$ Such a criterion essentially indicates that only the principal part influences the kernel's behavior, as adding lower order terms with bounded coefficients preserves any existing bounds. It is stable under $L^{\infty}$ perturbation of the principal part coefficients, in line with the second order case, and provides heat kernel estimates for operators with VMO coefficients. Based on this work and the strategy of Kim \cite{K}, we investigate the De Giorgi-Moser-Nash estimates and Gaussain estimates for higher-order parabolic operators in form of \begin{equation}\label{eq: ce00}\H:=\partial_{t} +\L :=\partial_{t} +\sum_{|\alpha|=|\beta|=m}(-1)^{|\alpha|}\partial^{\alpha}(a_{\alpha, \beta}(x)\partial^{\beta})\end{equation} in $\rr^{n+1}.$ 

A major goal of this paper is to prove that the fundamental solutions to parabolic operators shaped like \eqref{eq: ce00} with $t-$independent coefficients have Gaussian upper bounds and H\"{o}lder regularity estimates when $2m\geq n$ (Theorem \ref{theorem: ce1010}). We also show that if weak solutions of an elliptic equation ($\L$) verify the H\"{o}lder property $H(\mu, \infty)$ (the higher-order version of the De Giorgi-Moser-Nash estimates, see Definition \ref{definition: ce17}), then weak solutions of the corresponding parabolic equation 
($\H$) with $t-$independent coefficients satisfy a similar parabolic H\"{o}lder property (Theorem \ref{theorem: ce20}). Building upon this result, a small $L^{\infty}$ perturbation theorem (Theorem \ref{theorem: ce59}) is established for higher order parabolic equations with $t-$independent coefficients. This result accomplishes two important extensions: $(i)$ it generalizes the elliptic case established in \cite{AQ} to parabolic equations, and $(ii)$
it provides a higher-order analogue of the result in \cite{A}. It also yields that when the coefficients are $t-$independent and small $L^{\infty}$ perturbations of complex constants or VMO functions, weak solutions of a parabolic operator defined in \eqref{eq: ce00} verify the H\"{o}lder property $H(\mu, \infty)$ (see Remark \ref{remark: ce60}). We stress that these results can be extended to systems without difficulty. On the other hand, our main results have direct applications in developing higher-order layer potential theory, offering solutions to boundary value problems associated with higher-order parabolic operators. Readers may refer to \cite{ B1, B2, B3, NAS, NK} for more details on this topic.

The remaining sections are organized as follows. In Section 2 we introduce notations, basic definitions and assumptions, and function spaces. We then establish some technical lemmas for solutions to higher-order elliptic and parabolic equations with time-independent coefficients in Section 3, and use them to prove the H\"{o}lder property $H(\mu, \infty)$ in Section 4. The last section is devoted to building Gaussian estimates for the fundamental solution of parabolic equations under consideration.

\section{Preliminaries}
\noindent{\textbf{2.1. Notation.}} Given a typical point $X:=(t, x)\in \rdm:=\{(t, x)\in \rr^{n+1}: t\in \rr, x\in \rr^{n}\}$ with $x=(x_{1},...x_{n})$ ($n\geq 2$). For any $r>0,$ $B_{r}(x):=\{y\in \rz: |x-y|<r\}$ is an open ball in $\rz,$ and we use $Q_{r}:=Q_{r}(t, x)=(t-r^{2m}, t]\times B_{r}(x)$ and $Q^{*}_{r}(t, x):=[t, t+R^{2m})\times B_{R}(x)$ to denote half cubes in $\rdm.$ For any multi-index $\gamma=(\gamma_{1}, ...\gamma_{n}) \in \nn^{n},$ let $$D^{\gamma}:=\partial^{\gamma_{1}}_{x_{1}}... \partial^{\gamma_{n}}_{x_{n}}$$ and $|\gamma|:=\gamma_{1}+...+\gamma_{n}.$ We also let $D^{k}:=\nabla^{k}:=(D^{\gamma})_{|\gamma|=k}$ for any $0\leq k\leq m.$ Besides, we need the notation for averages, expressed as follows: $$\fint_{E}f:=\frac{1}{|E|}\int_{E}f(x)dx,$$ where $E$ is a Borel measurable set in $\rz$ or $\rdm$, with its Lebesgue measure $0<|E|<\infty,$ and $f$ is a locally integrable function.  

\noindent{\textbf{2.2. The parabolic operator of higher order.}} We consider homogeneous parabolic operators of order $2m$ ($m\geq 1$) defined as \eqref{eq: ce00}. The coefficients $\{a_{\alpha, \beta}\}$ are assumed to be bounded, uniformly elliptic (in the sense of G$\mathring{a}$rding inequality)  and independent of $t$, that is, for all $u\in H^{m}:=W^{m, 2}(\rz),$
 \begin{equation}\label{eq: ce200}\mbox{Re} <\L u, u>:=\sum_{|\alpha|=|\beta|=m}\int_{\rz}a_{\alpha, \beta}(x)\partial^{\alpha}u \overline{\partial^{\beta}u}\geq \frac{c_{1}}{2}\|\nabla^{m} u\|_{2}^{2}-\lambda_{0} \|u\|_{2}^{2},\end{equation} 

and  \begin{equation}\label{eq: ce1}|\sum_{|\alpha|=|\beta|=m}a_{\alpha, \beta}( x)\xi_{\alpha}\eta_{\beta}|\leq c_{2}\l(\sum_{|\beta|=m}|\xi_{\beta}|^{2}\r)^{1/2}\l(\sum_{|\beta|=m}|\xi_{\eta}|^{2}\r)^{1/2}, \quad \mbox{for all}\; x\in \rr^{n}.\end{equation} Below, in accordance with \cite{AQ}, we denote by $\Ez_{2m}(c_{1}, \lambda_{0}, c_{2})$ the class of homogeneous elliptic operators $\L$ of order $2m$ satisfying \eqref{eq: ce200} and \eqref{eq: ce1}. Moreover, we borrow the notation $\Ez_{2m}(c_{1}, \lambda_{0}, c_{2}, c_{3})$ from \cite[Definition 2]{AQ} to denote the inhomogeneous elliptic operators of order $2m.$ 

\noindent{\textbf{2.3. Weak solutions.}}  
Let $\Omega\subset \rz$ be a domain and $H^{m}(\Omega):=W^{m, 2}(\Omega)$ be the Sobolev space of complex valued functions $u,$ defined on $\Omega,$ such that $D^{\alpha} u\in L^{2}(\Omega, \cc)$ for any $|\alpha|\leq m.$ We say that $u\in W^{m, 2}(\Omega)$ is a weak solution of $$\L u=f \in L^{2}(\Omega) \;\quad \mbox{in}\; \Omega$$ if for any $\psi\in C_{0}^{\infty}(\Omega)$ $$\sum_{|\alpha|=|\beta|=m}\int_{\Omega}a_{\alpha, \beta}\partial^{\beta}u\overline{\partial^{\alpha}\psi}=\int_{\Omega} f\overline{\psi}.$$
Define $\Omega_{t_{1}, t_{2}}:=[t_{1}, t_{2}]\times \Omega \subset \rr\times \rz.$ 
We say that $u\in  H^{m, 2}(\Omega_{t_{1}, t_{2}}):=L^{2}([t_{1}, t_{2}], W^{m, 2}(\Omega))$ is a weak solution of $\H$ in $\Omega_{t_{1}, t_{2}},$ if for any $\phi\in C_{0}^{\infty}(\Omega_{t_{1}, t_{2}}),$ $$-\int_{\Omega_{t_{1}, t_{2}}} u\overline{\partial_{t}\phi}+\sum_{|\alpha|=|\beta|=m}\int_{\Omega_{t_{1}, t_{2}}}a_{\alpha, \beta}\partial^{\beta}u\overline{\partial^{\alpha}\phi}=0.$$ By the definition of the weak solution to $\H,$ the adjoint $\H^{*}$ of $\H$ can be identified formally as the backward-in-time operator $$-\partial_{t}+\L^{*}:=-\partial_{t}+\sum_{|\alpha|=|\beta|=m}(-1)^{|\beta|}\partial^{\beta}(\overline{a_{\alpha, \beta}}\partial^{\alpha}).$$ Clearly, if $\L\in \Ez_{2m}(c_{1}, 0, c_{2}),$ so is $\L^{*}.$
To circumvent technical difficulties, we impose a smoothness assumption on the coefficients-a widely used strategy in numerous studies such as \cite{DK, HK, K}. Thus our weak solutions are indeed strong (classical) solutions. 
This simplification is ultimately harmless, as our arguments never quantitatively rely on the additional regularity. This is why we make the dependence of constants clear in what follows. 

\noindent{\textbf{2.4. Basic definitions.}} According to \cite[Theorem 3.10]{B1} and \cite[Definition 6]{AQ} (the global case), we can introduce the following definition, that is the higher-order version of the De Giorgi-Moser-Nash estimates in \cite{AAA, NAS, HMM1, NK}.
  
\begin{definition}\label{definition: ce17} Suppose $\L\in \Ez_{2m}(c_{1}, 0, c_{2}).$ Let $\mu \in (\max\{0, m-\frac{n}{2}\}, m]\setminus \nn.$ Write $\mu=l+\nu,$ where $l\in \{0,1,...m-1\}$ and $\nu\in (0, 1).$ We say that $\L$ has the H\"{o}lder property $H(\mu, \infty)$ if the following conditions hold: there exists a constant $c_{0}>0$ such that for all $0<R<\infty,$ for all $x_{0}\in \rz,$ for all solutions $u$ of $\L$ in $B_{R}(x_{0})$ and for all multi-index $\gamma \in \nn^{n}$ with $|\gamma|\leq l,$
\begin{equation}\label{eq: ce18}
R^{|\gamma|}\sup_{B_{R/2}(x_{0})}|D^{\gamma} u|\leq c_{0} \l(\fint_{B_{R}(x_{0})}|u|^{2}\r)^{1/2},
\end{equation} and if $|\gamma|=l,$
\begin{equation}\label{eq: ce19}
R^{|\gamma|+\nu}[[D^{\gamma} u]]_{\nu, B_{R/2}(x_{0})}\leq c_{0} \l(\fint_{B_{R}(x_{0})}|u|^{2}\r)^{1/2},
\end{equation} where $$[[u]]_{\nu, \Omega}:=\sup_{x\neq y, x, y\in \Omega}\frac{|u(x)-u(y)|}{|x-y|^{\nu}}.$$
Similarly, letting $Q_{R}:=Q_{R}(t_{0}, x_{0}),$ we say that $\H$ defined in \eqref{eq: ce00} satisfies the H\"{o}lder property $H(\mu, \infty)$ if \eqref{eq: ce18}-\eqref{eq: ce19} hold with $B_{R/2}(x_{0}), B_{R}(x_{0})$ and $[[D^{\gamma} u]]_{\nu, B_{R/2}(x_{0})}$ replaced by $Q_{R/2}, Q_{R}$ and $[D^{\gamma} u]_{\nu, Q_{R/2}},$ respectively, where $$[u]_{\nu, Q_{R}}:=\sup_{X\neq Y, X, Y\in Q_{R}}\frac{|u(X)-u(Y)|}{|X-Y|^{\nu}},\quad |X-Y|:=|t-s|^{\frac{1}{2m}}+ |x-y|.$$
Moreover, we say that $\H^{*}$ satisfies the H\"{o}lder property $H(\mu, \infty)$ if \eqref{eq: ce18}-\eqref{eq: ce19} hold with $B_{R/2}(x_{0})$ $ B_{R}(x_{0})$ and $[[D^{\gamma} u]]_{\nu, B_{R/2}(x_{0})}$ replaced by $Q^{*}_{R/2}, Q^{*}_{R}$ and $[D^{\gamma} u]_{\nu, Q^{*}_{R/2}},$ respectively. 
\end{definition}

We also introduce the so-called Gaussain estimates, that is the global property defined in \cite[Definition 9]{AQ}.

\begin{definition}\label{definition: ce3001} Assume $\L \in \Ez_{2m}(c_{1}, \lambda_{0}, c_{2}, c_{3}).$ Let $\mu \in (\max\{0, m-\frac{n}{2}\}, m]\setminus \nn.$ Write $\mu=l+\nu,$ where $l\in \{0,1,...m-1\}$ and $\nu\in (0, 1).$ We say that $\L$ has the Gaussain property $G(\mu, \infty)$ if the following condition holds: there exist constant $c>0$ and $a>0$ such that for all $0<t<\infty,$ for all $x, y, h\in \rz,$ and for all multi-index $\gamma \in \nn^{n}$ with $|\gamma|\leq l,$
\begin{equation}\label{eq: ce3002}
|D_{x}^{\gamma}K_{t}(x, y)|+|D_{y}^{\gamma}K_{t}(x, y)|\leq \frac{c}{(t-s)^{\frac{n+|\gamma|}{2m}}}g_{m, a}\l(\frac{|x-y|}{t^{\frac{1}{2m}}}\r),
\end{equation} and if $|\gamma|=l,$
\begin{equation}\label{eq: ce3003}
|D_{x}^{\gamma}K_{t}(x+h, y)-D_{x}^{\gamma}K_{t}(x, y)|\leq \frac{c}{t^{\frac{n+l}{2m}}}\l(\frac{|h|}{\frac{1}{2m}}\r)^{\nu},
\end{equation} 
\begin{equation}\label{eq: ce3004}
|D_{y}^{\gamma}K_{t}(x, y+h)-D_{y}^{\gamma}K_{t}(x, y)|\leq \frac{c}{t^{\frac{n+l}{2m}}}\l(\frac{|h|}{\frac{1}{2m}}\r)^{\nu},
\end{equation} where $g_{m, a}(s):=e^{-as^{\frac{2m}{2m-1}}}$ for $s>0$ and $K_{t}(, y)$ denotes the distributional kernel of the semigroup $e^{-t\L}.$
\end{definition}

\noindent{\textbf{2.5. The explicit constants.}}  The two letters $C, c$ will be used to denote universal constants depending on $n, m, c_{1}, c_{2}$ and they may vary from line to line. To emphasize the extra dependence on the prescribed quantities $a, b,....,$ we write $C=C(a, b,....)$ and $c=c(a, b,....).$

\section{Technical lemmas}

This section focuses on constructing the key lemmas that will serve as the foundation for the proof of the H\"{o}lder property from elliptic equations to parabolic equations in the subsequent section. We always assume $\L \in \Ez_{2m}(c_{1}, 0, c_{2})$ in this section.  

\begin{lemma}\label{lemma: ce3} Let $R>0$ and $u\in H^{m, 2}(Q_{R}(t, x))$ be a weak solution of $\H$ in $Q_{R}(t, x).$ Then for any $0<r<R,$ we have 
\begin{equation}\label{eq: ce4}\sup_{t-r^{2m}<s\leq t}\|u(s, \cdot)\|^{2}_{L^{2}(B_{r}(x))}+\int_{Q_{r}(t, x)}|\nabla^{m} u|^{2}\leq C \frac{\int_{Q_{R}(t, x)}|u|^{2}}{(R-r)^{2m}}.\end{equation} In particular, for $j=0,1, 2,...m-1,$ \begin{equation}\label{eq: ce5}\int_{Q_{r}(t, x)}|\nabla^{j} u|^{2}\leq C \frac{\int_{Q_{R}(t, x)}|u|^{2}}{(R-r)^{2j}}.\end{equation} 
\end{lemma}
{\it Proof.}\quad For simplicity, we assume $t=0, x=0.$ Set $$r_{0}=r, r_{k}=r+\sum_{i=1}^{k}\frac{R-r}{2^{k}},\;\mbox{k=1,2,...}$$ and $$s_{k}=\frac{r_{k}+r_{k+1}}{2},\;\mbox{k=0,1,2,...}.$$ We choose a sequence of nonnegative functions $\phi_{k}\in C_{0}^{\infty}(\rdm)$ such that $\phi_{k}=1$ on $Q_{r_{k}}$ and $\phi_{k}=0$ on $\rdm\setminus (-s_{k}^{2m}, s_{k}^{2m})\times B_{s_{k}},$ moreover, $$\|\partial_{t}\phi_{k}\|_{L^{\infty}}\leq c\frac{2^{2mk}}{(R-r)^{2m}}, \quad \|D^{l}\phi_{k}\|_{L^{\infty}}\leq c \frac{2^{2lk}}{(R-r)^{l}}, \;\mbox{l=0, 1,2,...}.$$ Testing the equation by $u\phi_{k}^{2}1_{(-\infty, s)}(t)$ for any $-r^{2m}<s\leq 0$ and using the fact that the coefficients are $t-$independent, we get \begin{equation}\label{eq: ce6}-\int_{Q_{R}\cap E_{s}}|u|^{2}\phi_{k}\partial_{t}\phi_{k}+\frac{1}{2}\int_{B_{R}(0)}|u\phi_{k}|^{2}(s, x)dx+\sum_{|\alpha|=|\beta|=m}\int_{Q_{R}\cap E_{s}}a_{\alpha, \beta}\partial^{\beta}u\partial^{\alpha}(\overline{u}\phi_{k}^{2})=0,\end{equation} where $E_{s}:=(-\infty, s).$ By Leibniz's rule, the equation \eqref{eq: ce6} can be rewritten as  \begin{equation}\label{eq: ce7}\begin{aligned}
\int_{Q_{R}\cap E_{s}}a_{\alpha, \beta}\partial^{\beta}u\partial^{\alpha}(\overline{u}\phi_{k}^{2})&=\int_{Q_{R}\cap E_{s}}a_{\alpha, \beta}\partial^{\beta}u\partial^{\alpha}(\overline{u}\phi_{k}^{2})-\int_{Q_{R}\cap E_{s}}a_{\alpha, \beta}\partial^{\beta}(u\phi_{k})\overline{\partial^{\alpha} (u\phi_{k})}\\[4pt]
&\quad +\int_{Q_{R}\cap E_{s}}a_{\alpha, \beta}\partial^{\beta}(u\phi_{k})\overline{\partial^{\alpha} (u\phi_{k})}:=I_{1}+I_{2}.\end{aligned}
\end{equation} Bear in mind that $\phi(s, \cdot)\in C_{0}^{\infty}(B_{\tau}(x_{0}))$ for all $s\leq 0.$ Then, it follows from \eqref{eq: ce200} that \begin{equation}\label{eq: ce8}\mbox{Re}\; I_{2}\geq c_{1}\int_{Q_{R}\cap E_{s}} | \nabla^{m}( u\phi_{k} )|^{2}.\end{equation} To proceed, we decompose $I_{1}$ into 
\begin{equation*}\begin{aligned}I_{1}&=\int_{Q_{R}\cap E_{s}}a_{\alpha, \beta}[\phi_{k}\partial^{\beta}u -\partial^{\beta}(u\phi_{k})]\overline{\partial^{\alpha}(u\phi_{k})}+a_{\alpha, \beta}[\phi_{k}\overline{\partial^{\alpha}(u\phi_{k})} -\overline{\partial^{\alpha}(u\phi_{k}^{2})}]\partial^{\beta}u\\[4pt]
&=\sum_{\gamma<\beta}C_{\beta}^{\gamma}\int_{Q_{R}\cap E_{s}}a_{\alpha, \beta}\partial^{\gamma}u \partial^{\beta-\gamma}\phi_{k}\overline{\partial^{\alpha}(u\phi_{k})}-\sum_{\tau< \alpha}C_{\alpha}^{\tau}\partial^{\beta}u\partial^{\alpha-\tau}\phi_{k}\overline{\partial^{\tau}(u\phi_{k})}\\[4pt]
&:=I_{11}-I_{12}.
\end{aligned}
\end{equation*} Recalling the definition of $\phi_{k},$ the term $I_{11}$ can be further written as \begin{equation*}\begin{aligned}I_{11}
&=\sum_{\gamma<\beta}C_{\beta}^{\gamma}\int_{Q_{R}\cap E_{s}}a_{\alpha, \beta}\partial^{\gamma}(u\phi_{k+1}) \partial^{\beta-\gamma}\phi_{k}\overline{\partial^{\alpha}(u\phi_{k+1}\phi_{k})}\\[4pt]
&=\sum_{\gamma<\beta}\sum_{\tau\leq \alpha}C_{\alpha}^{\tau}C_{\beta}^{\gamma} \int_{Q_{R}\cap E_{s}}a_{\alpha, \beta}\partial^{\gamma}(u\phi_{k+1}) \partial^{\beta-\gamma}\phi_{k}\partial^{\alpha-\tau}\phi_{k}\overline{\partial^{\tau}(u\phi_{k+1})}.
\end{aligned}
\end{equation*} 
In veiw of the property of $\phi_{k}$ and interpolation inequalities, we can estimate $I_{11}$ as follows: 

\begin{equation*}\begin{aligned}
|I_{111}|&\leq  c\sum_{\gamma<\beta}\sum_{\tau\leq \alpha}C_{\alpha}^{\tau}C_{\beta}^{\gamma} \l(\frac{2^{k}}{R-r}\r)^{m-|\gamma|+m-|\tau|}\|\partial^{\gamma}(u\phi_{k+1})\|_{L^{2}(R^{n+1}_{0})}\|\partial^{\tau}(u\phi_{k+1})\|_{L^{2}(R^{n+1}_{0})}\\[4pt]
\end{aligned}
\end{equation*}
\begin{equation*}\begin{aligned}
&\leq  c\sum_{\gamma<\beta}\sum_{\tau\leq \alpha}C_{\alpha}^{\tau}C_{\beta}^{\gamma} \l(\frac{2^{k}}{R-r}\r)^{m-|\gamma|}\|D^{m}(u\phi_{k+1})\|_{L^{2}(R^{n+1}_{0})}^{|\gamma|/m}\|u\phi_{k+1}\|_{L^{2}(R^{n+1}_{0})}^{1-\frac{|\gamma|}{m}}\\[4pt]
&\quad\quad\times \l(\frac{2^{k}}{R-r}\r)^{m-|\tau|}\|D^{m}(u\phi_{k+1})\|_{L^{2}(R^{n+1}_{0})}^{|\tau|/m}\|u\phi_{k+1}\|_{L^{2}(R^{n+1}_{0})}^{1-\frac{|\tau|}{m}}\\[4pt]
&\leq  c\sum_{\gamma<\beta}\sum_{\tau\leq \alpha}C_{\alpha}^{\tau}C_{\beta}^{\gamma}\l(\ez\|D^{m}(u\phi_{k+1})\|_{L^{2}(R^{n+1}_{0})}+\l(\frac{2^{k}}{R-r}\r)^{m}\ez^{\frac{|\gamma|}{|\gamma|-m}}\|u\phi_{k+1}\|_{L^{2}(R^{n+1}_{0})}\r)\\[4pt]
&\quad\quad\times \l(\ez\|D^{m}(u\phi_{k+1})\|_{L^{2}(R^{n+1}_{0})}+\l(\frac{2^{k}}{R-r}\r)^{m}\ez^{\frac{|\tau|}{|\tau|-m}}\|u\phi_{k+1}\|_{L^{2}(R^{n+1}_{0})}\r)
\\[4pt]
&\leq \ez \|D^{m}(u\phi_{k+1})\|_{L^{2}(R^{n+1}_{0})}^{2}+c\l(\frac{2^{k}}{R-r}\r)^{2m}\sum_{l=0}^{m-1}\ez^{\frac{l}{l-m}}\|u\phi_{k+1}\|_{L^{2}(R^{n+1}_{0})}^{2},
\end{aligned}
\end{equation*}
where $R^{n+1}_{0}:=(-\infty, 0)\times \rz.$ We remark that, in the above argument, we ignore the case $|\tau|=m,$ which, in fact, is easier. Observing that $I_{12}$ and $I_{11}$ are of the same type, thus, by the way to tackle $I_{11}$ we conclude  $$I_{12}\leq \ez \|D^{m}(u\phi_{k+1})\|_{L^{2}(R^{n+1}_{0})}^{2}+c\l(\frac{2^{k}}{R-r}\r)^{2m}\sum_{l=0}^{m-1}\ez^{\frac{l}{l-m}}\|u\phi_{k+1}\|_{L^{2}(R^{n+1}_{0})}^{2}:=I_{3}.$$
Clearly, the same bound holds for $I_{1}.$

Collecting \eqref{eq: ce6}-\eqref{eq: ce8} and the above bound for $I_{1},$ we obtain 
\begin{equation}\label{eq: ce9}\begin{aligned}
\frac{1}{2}\int_{B_{R}(0)}|u\phi_{k}|^{2}(s, x)dx+c_{1}\int_{Q_{R}\cap E_{s}} |\nabla^{m} (u\phi_{k})|^{2}&\leq \int_{Q_{R}\cap E_{s}}|u|^{2}|\phi_{k}\partial_{t}\phi_{k}| +c I_{3}\\[4pt]
&\leq  c \frac{2^{2mk}}{(R-r)^{2m}}\|u\|_{L^{2}(Q_{R})}^{2}+cI_{3}.
\end{aligned}
\end{equation}  
Letting $s=0$ in \eqref{eq: ce9} we have \begin{equation*}
\frac{1}{2}\sup_{-r^{2m}<s\leq 0}\int_{B_{R}(0)}|u\phi_{k}|^{2}(s, x)dx+\int_{R^{n+1}_{0}} |\nabla^{m} (u \phi_{k} )|^{2}\leq c \frac{2^{2mk}}{(R-r)^{2m}}\|u\|_{L^{2}(Q_{R})}^{2} +c I_{3},
\end{equation*}
which implies  \begin{equation}\label{eq: ce11}\begin{aligned}
\sup_{-r^{2m}<s\leq 0}\int_{B_{R}(0)}|u\phi_{k}|^{2}(s, x)dx&+\int_{R^{n+1}_{0}} | \nabla^{m} (u \phi_{k})|^{2}\\[4pt]
&\leq  \ez \int_{R^{n+1}_{0}}|\nabla^{m}(u\phi_{k+1})|^{2}+c\sum_{l=0}^{m-1}(1+\ez^{\frac{l}{l-m}})\frac{2^{2mk}}{(R-r)^{2m}}\|u\|_{L^{2}(Q_{R})}^{2}.
\end{aligned}
\end{equation}
Note that $\phi_{0}=1$ on $Q_{r}.$ If we let $k=0$ in \eqref{eq: ce11}, then \begin{equation}\label{eq: ce12}\begin{aligned}
\sup_{-r^{2m}<s\leq 0}\int_{B_{r}(0)}|u|^{2}(s, x)dx
&\leq c  \int_{R^{n+1}_{0}}|\nabla^{m}(u\phi_{1})|^{2}+\frac{c}{(R-r)^{2m}}\|u\|_{L^{2}(Q_{R})}^{2}.
\end{aligned}
\end{equation} On the other hand, due to \eqref{eq: ce11} again, an iteration argument displayed in \cite[Lemma 3.2]{DK} leads to \begin{equation}\label{eq: ce5010}\int_{R^{n+1}_{0}}|\nabla^{m}(u\phi_{0})|^{2}+\int_{R^{n+1}_{0}}|\nabla^{m}(u\phi_{1})|^{2}\leq  \frac{c}{(R-r)^{2m}}\|u\|_{L^{2}(Q_{R})}^{2}.\end{equation} Thus, by \eqref{eq: ce12}, we arrive at $$\sup_{-r^{2m}<s\leq 0}\int_{B_{r}(0)}|u|^{2}(s, x)dx\leq \frac{c}{(R-r)^{2m}}\|u\|_{L^{2}(Q_{R})}^{2}.$$ In particular, by \eqref{eq: ce5010} as well as interpolation inequalities, \eqref{eq: ce5} is almost trivial. Consequently, we complete the proof for \eqref{eq: ce4} and \eqref{eq: ce5}.

\hfill$\Box$

\begin{lemma}\label{lemma: ce21} Let $u\in H^{m, 2}(Q_{2R})$ be a solution of $\H$ in $Q_{2R}:=Q_{2R}(t_{0}, x_{0}),$ then $u_{t}\in L^{2}(Q_{r})$ for any $0<r<R,$ and the estimate 
\begin{equation}\label{eq: ce22}\|u_{t}\|_{L^{2}(Q_{r})}\leq C \frac{\|\nabla ^{m} u\|_{L^{2}(Q_{2R})}}{(R-r)^{m}}\end{equation} holds. Moreover, if $u$ is a solution of $\H$ in $Q_{2r},$ then \eqref{eq: ce22} along with the energy estimate yield \begin{equation}\label{eq: ce23}\|u_{t}\|_{L^{2}(Q_{r})}\leq C \frac{\|u\|_{L^{2}(Q_{2r})}}{r^{2m}}.\end{equation} 
\end{lemma}

{\it Proof.}\quad Given two positive constants $\sigma, \tau$ such that $\sigma<\tau\leq R.$ Let $\phi$ be a smooth cut-off function such that $\phi\equiv 1$ in $Q_{\sigma},$ vanishing near the boundary and verifying \begin{equation}\label{eq: ce29}0\leq \phi\leq 1, \quad \|\phi_{t}\|_{L^{\infty}}\leq \frac{c}{(\tau-\sigma)^{2m}} \quad\mbox{and}\quad \|D^{\gamma}\phi\|_{L^{\infty}}\leq \frac{c}{(\tau-\sigma)^{|\gamma|}}, \; \forall\; |\gamma|\leq m.\end{equation} For simplicity, we set $Q_{\tau, s}:={s}\times B_{\tau}(x_{0})$ with $s\in (t_{0}-\tau^{2m}, t_{0}).$ Of course, on each slice $Q_{\tau, s}$,  
\begin{equation*}\begin{aligned}
0=\int_{Q_{\tau, s}}\phi^{2}|u_{t}|^{2}+a_{\alpha, \beta}\partial^{\beta}u\phi\overline{\partial^{\alpha} (u_{t} \phi)}+\int_{Q_{\tau, s}} \sum_{\alpha_{1}+\alpha_{2}=\alpha, |\alpha_{2}|<m}a_{\alpha, \beta}\partial^{\beta}u(D^{\alpha_{1}} \phi)\overline{D^{\alpha_{2}}(u_{t} \phi)}
\end{aligned}
\end{equation*} holds. 
From this equality and Young's inequality, it follows that 
\begin{equation}\label{eq: ce61}\begin{aligned}
\int_{Q_{\tau, s}}\phi^{2}|u_{t}|^{2}&\leq \frac{\ez}{2}\int_{Q_{\tau, s}}  |\nabla^{m} (u_{t} \phi)|^{2}+\frac{c}{\ez}\int_{Q_{\tau, s}} |\nabla^{m}u\phi|^{2}\\[4pt]
&+\int_{Q_{\tau, s}} \sum_{\alpha_{1}+\alpha_{2}=\alpha, |\alpha_{2}|<m}|D^{\alpha_{1}} \phi| |\partial^{\beta}u||D^{\alpha_{2}}(u_{t}\phi)|:=J_{1}+J_{2}+J_{3}.
\end{aligned}
\end{equation}
We first handle $J_{3}.$ Using Young's inequality again and interpolation inequalities, it is easy to deduce \begin{equation*}\begin{aligned}J_{3}&\leq \frac{c}{\ez}\int_{Q_{\tau, s}} |\nabla^{m}u|^{2}+c\ez\int_{Q_{\tau, s}} \sum_{\alpha_{1}+\alpha_{2}=\alpha, |\alpha_{2}|<m}|D^{\alpha_{1}} \phi|^{2}|D^{\alpha_{2}}(u_{t}\phi)|^{2}\\[4pt]
&\leq \frac{c}{\ez}\int_{Q_{\tau, s}} |\nabla^{m}u|^{2}+ c\frac{\ez }{2}\int_{Q_{\tau, s}}  |\nabla^{m} (u_{t} \phi)|^{2}+c\frac{\ez}{(\tau-\sigma)^{2m}}\int_{Q_{\tau, s}} |u_{t}\phi|^{2}.\end{aligned}
\end{equation*} To manage $J_{1}$, we can apply the energy estimates \eqref{eq: ce4}-\eqref{eq: ce5} (adjusting the support of $\phi$ if necessary) to derive  
\begin{equation}\label{eq: ce24}\int_{Q_{\tau}}  |\nabla^{m} (u_{t} \phi)|^{2}\leq \frac{c}{(\tau-\sigma)^{2m}}\int_{Q_{\tau}} |u_{t}|^{2}\end{equation} since $u_{t}$ is also a weak solution of $\H$ in $Q_{2R}.$ $J_{2}$ is simply arranged by diacarding $\phi.$ Consequently, \eqref{eq: ce61} contributes to\begin{equation*}\begin{aligned}
\int_{Q_{\sigma}}|u_{t}|^{2}&\leq \frac{\ez}{2}\frac{c}{(\tau-\sigma)^{2m}}\int_{Q_{\tau}} |u_{t}|^{2}+\frac{c}{\ez}\int_{Q_{\tau}} |\nabla^{m}u|^{2}\\[4pt]&+ \frac{\ez }{2}\frac{c}{(\tau-\sigma)^{2m}}\int_{Q_{\tau}} |u_{t}|^{2}+c\frac{\ez}{(\tau-\sigma)^{2m}}\int_{Q_{\tau}} |u_{t}|^{2}.
\end{aligned}
\end{equation*} Letting $\ez:=(\tau-\sigma)^{2m}\ez_{1}$ with $\ez_{1}=\ez_{1}(c_{1}, c_{2}, m, n)$ small enough in the latter inequality we have 

\begin{equation}\label{eq: ce25} \int_{Q_{\sigma}}|u_{t}|^{2}\leq \frac{1}{2} \int_{Q_{\tau}} |u_{t}|^{2}+\frac{c}{(\tau-\sigma)^{2m}}\int_{Q_{\tau}} |\nabla^{m}u|^{2}.\end{equation} Armed with \eqref{eq: ce25}, an application of the iteration lemma \cite[Lemma 5.1]{GM} yields\begin{equation}\label{eq: ce30}\int_{Q_{r}}|u_{t}|^{2}\leq \frac{c}{(R-r)^{2m}}\int_{Q_{R}} |\nabla^{m}u|^{2}\quad \mbox{for all}\; 0<r<R.\end{equation} The proof is complete.

\hfill$\Box$

As a result of Lemma \ref{lemma: ce21}, we can deduce the following uniform estimates on the slices.

\begin{lemma}\label{lemma: ce26} Let $u\in H^{m, 2}(Q_{2r})$ be a solution of $\H$ in $Q_{2r}:=Q_{2r}(t_{0}, x_{0}).$ Then $\nabla^{j}u(s, \cdot), u_{t}(s, \cdot)\in L^{2}(B_{r})$ for all $0\leq j\leq m$ and $ s\in [t_{0}-r^{2m}, t_{0}],$ and satisfy 
\begin{equation}\label{eq: ce27}\sup_{s\in [t_{0}-r^{2m}, t_{0}]}\|\nabla^{j}u(s, \cdot)\|_{L^{2}(B_{r})}\leq C \frac{\| u\|_{L^{2}(Q_{2r})}}{r^{m+j}},\end{equation}\begin{equation}\label{eq: ce28}\sup_{s\in [t_{0}-r^{2m}, t_{0}]}\|u_{t}(s, \cdot)\|_{L^{2}(B_{r})}\leq C \frac{\| u\|_{L^{2}(Q_{2r})}}{r^{3m}}.\end{equation} 
\end{lemma}

{\it Proof.}\quad Fix $r\leq \sigma<\tau\leq 2r$ and let $\phi$ be a smooth cut-off function defined as in \eqref{eq: ce29}. We first note that $$\sup_{t-\sigma^{2m}<s\leq t}\|u_{t}(s, \cdot)\|^{2}_{L^{2}(B_{\sigma}(x))}\leq c \frac{\int_{Q_{\sigma+\frac{\tau -\sigma}{2}}}|u_{t}|^{2}}{(\tau-\sigma)^{2m}}$$ holds thanks to the $t$-independence of the coefficients and the energy estimate \eqref{eq: ce4}. Furthermore, it follows from \eqref{eq: ce30} that $$\int_{Q_{\sigma+\frac{\tau -\sigma}{2}}}|u_{t}|^{2} \leq c \frac{\|u\|^{2}_{L^{2}(Q_{\tau})}}{(\tau-\sigma)^{4m}}.$$ The two estimates imply \begin{equation}\label{eq: ce31}\sup_{t-\sigma^{2m}<s\leq t}\|u_{t}(s, \cdot)\|^{2}_{L^{2}(B_{\sigma}(x))}\leq c \frac{\|u\|^{2}_{L^{2}(Q_{\tau})}}{(\tau-\sigma)^{6m}}.\end{equation}

 For the moment, we assume that $u$ is a weak solution of $\H$ in $Q_{4r}.$ Observe that, for any fixed $s,$ \begin{equation*}\begin{aligned}0&=\int_{B_{\tau}}\phi^{2}u_{t}(s, \cdot)\overline{u(s, \cdot)}+a_{\alpha, \beta}\partial^{\beta}(u\phi)\overline{\partial^{\alpha} (u\phi) }+\int_{B_{\tau}}a_{\alpha, \beta}\partial^{\beta}u[\phi\overline{D^{\alpha}(u\phi) }-\overline{D^{\alpha}(u\phi^{2}) }]\\[4pt]
 &\quad\quad+\int_{B_{\tau}}a_{\alpha, \beta}[\phi\partial^{\beta}u -\partial^{\beta}(u\phi)]\overline{\partial^{\alpha}(u\phi)}\\[4pt]
\end{aligned}
 \end{equation*}  
\begin{equation*}\begin{aligned} 
&=\sum_{\gamma<\beta}C_{\beta}^{\gamma}\int_{B_{\tau}}a_{\alpha, \beta}\partial^{\gamma}u \partial^{\beta-\gamma}\phi\overline{\partial^{\alpha}(u\phi)}-\sum_{\tau< \alpha}C_{\alpha}^{\tau}\int_{B_{\tau}}a_{\alpha, \beta}\partial^{\beta}u\partial^{\alpha-\tau}\phi\overline{\partial^{\tau}(u\phi)}\\[4pt]
&\quad\quad +\int_{B_{\tau}}\phi^{2}u_{t}(s, \cdot)\overline{u(s, \cdot)}+a_{\alpha, \beta}\partial^{\beta}(u\phi)\overline{\partial^{\alpha} (u\phi) }.  \end{aligned}
 \end{equation*} By this equality and \eqref{eq: ce200}, we see \begin{equation}\label{eq: ce3010}\begin{aligned}c_{1}\int_{B_{\tau}} |\nabla^{m} (u\phi)|^{2}&\leq c \int_{B_{\tau}}\phi^{2}|u_{t}(s, \cdot)||u(s, \cdot)|+c\sum_{\gamma<\beta}C_{\beta}^{\gamma}\int_{B_{\tau}}|\partial^{\gamma}u| |\partial^{\beta-\gamma}\phi||\partial^{\alpha}(u\phi)|\\[4pt]
 \end{aligned}
\end{equation}
\begin{equation*}\begin{aligned} 
&\quad\quad + c\sum_{\tau< \alpha}C_{\alpha}^{\tau}\int_{B_{\tau}}|\partial^{\beta}u| |\partial^{\alpha-\tau}\phi| |\partial^{\tau}(u\phi)|:=M_{1}+M_{2}+M_{3}.\end{aligned}
\end{equation*}
Exploiting in succession Young's inequality, \eqref{eq: ce4} and \eqref{eq: ce31} we attain  \begin{equation}\label{eq: ce3011}\begin{aligned}M_{1}&\lesssim \ez \int_{B_{\tau}}|u_{t}(s, \cdot)|^{2}\phi^{2}+\frac{1}{\ez}\int_{B_{\tau}}|u(s, \cdot)|^{2}\phi^{2}\\[4pt]
&\leq \frac{\ez }{(\tau-\sigma)^{6m}}\int_{Q_{4r}}|u|^{2}+\frac{c}{\ez(\tau-\sigma)^{2m}}\int_{Q_{4r}}|u|^{2}.\end{aligned}
\end{equation} For the term $M_{2}$, by interpolation inequalities and \eqref{eq: ce4}, we have \begin{equation*}\begin{aligned}M_{2}&\leq \frac{c_{1}}{8} \int_{B_{\tau}}|\nabla^{m}((u\phi)(s, \cdot) )|^{2}+c\ez_{1}\int_{B_{\tau}}|\nabla^{m}u(s, \cdot)|^{2}+\frac{c(\ez_{1})}{(\tau-\sigma)^{2m}}\int_{B_{\tau}\cap \supp \phi}|u(s, \cdot)|^{2}\\[4pt]
&\leq \frac{c_{1}}{8} \int_{B_{\tau}}|\nabla^{m}((u\phi)(s, \cdot))|^{2}\phi^{2}+c\ez_{1}\int_{B_{\tau}}|\nabla^{m}u(s, \cdot)|^{2}+\frac{c(\ez_{1})}{(\tau-\sigma)^{4m}}\int_{Q_{4r}}|u|^{2}.\end{aligned}
\end{equation*}  By the same token, \begin{equation*}\begin{aligned}M_{3}&\leq \frac{c_{1}}{4} \int_{B_{\tau}}|\nabla^{m}u(s, \cdot)|^{2}+c\ez_{2}\int_{B_{\tau}}|\nabla^{m}((u\phi)(s, \cdot))|^{2}+\frac{c(\ez_{2})}{(\tau-\sigma)^{2m}}\int_{B_{\tau}\cap \supp \phi}|u(s, \cdot)|^{2}\\[4pt]
&\leq \frac{c_{1}}{4} \int_{B_{\tau}}|\nabla^{m}((u\phi)(s, \cdot))|^{2}+c\ez_{2}\int_{B_{\tau}}|\nabla^{m}u(s, \cdot)|^{2}+\frac{c(\ez_{2})}{(\tau-\sigma)^{4m}}\int_{Q_{4r}}|u|^{2}.\end{aligned}
\end{equation*}
Summarizing the above inequalities and letting $\ez:=(\tau-\sigma)^{2m}\ez_{1}$ in \eqref{eq: ce3011}, we derive from \eqref{eq: ce3010} that 
\begin{equation}\label{eq: ce70}\begin{aligned}\int_{B_{\tau}} |\nabla^{m} (u\phi)|^{2}\leq c \ez_{1} \int_{B_{\tau}}|\nabla^{m}u(s, \cdot)|^{2}+c\ez_{2}\int_{B_{\tau}}|\nabla^{m}u(s, \cdot)|^{2}+\frac{c(\ez_{1}, \ez_{2})}{(\tau-\sigma)^{4m}}\int_{Q_{4r}}|u|^{2}.\end{aligned}
\end{equation} In \eqref{eq: ce70}, by taking the supremum with respect to $t-\sigma^{2m}<s\leq t$ and choosing $\ez_{1}, \ez_{2}, $ small enough, we obtain $$\sup_{t-\sigma^{2m}<s\leq t}\|\nabla^{m}u(s, \cdot)\|^{2}_{L^{2}(B_{\sigma}(x))}\leq \theta \sup_{t-\tau^{2m}<s\leq t}\|\nabla^{m}u(s, \cdot)\|^{2}_{L^{2}(B_{\tau}(x))}+\frac{c}{(\tau-\sigma)^{4m}}\int_{Q_{4r}}|u|^{2},$$ where $r\leq \sigma<\tau\leq 2r$ and $0<\theta<1.$ As a consequence, an application of the iteration lemma \cite[Lemma 5.1]{GM} leads to \begin{equation}\label{eq: ce32}\sup_{t-r^{2m}<s\leq t}\|\nabla^{m}u(s, \cdot)\|^{2}_{L^{2}(B_{r}(x))}\leq c \frac{\|u\|^{2}_{L^{2}(Q_{4r})}}{r^{4m}}.\end{equation} Then, a well-known covering argument yields \eqref{eq: ce27} when $j=m.$ Recall that \eqref{eq: ce4} implies \begin{equation}\label{eq: ce33}\sup_{t-r^{2m}<s\leq t}\|u(s, \cdot)\|^{2}_{L^{2}(B_{r}(x))}\lesssim \frac{\|u\|^{2}_{L^{2}(Q_{2r})}}{r^{2m}}.\end{equation} Connecting \eqref{eq: ce32} and \eqref{eq: ce33}, also using interpolation inequalities (\cite[Theorem 2.1.14]{JK}), we eventually arrive at \begin{equation*}\begin{aligned}\sup_{t-r^{2m}<s\leq t}\|\nabla^{j}u(s, \cdot)\|^{2}_{L^{2}(B_{r}(x))}&\leq c r^{-2j}\sup_{t-r^{2m}<s\leq t}\|u(s, \cdot)\|^{2}_{L^{2}(B_{r}(x))}\\[4pt]
&\quad\quad+c\sup_{t-r^{2m}<s\leq t}r^{2(m-j)}\|\nabla^{m}u(s, \cdot)\|^{2}_{L^{2}(B_{r}(x))}\\[4pt]
&\leq c \frac{\|u\|^{2}_{L^{2}(Q_{2r})} } {r^{2m+2j}}.
 \end{aligned}
\end{equation*}
The lemma is proved. 

\hfill$\Box$

The following lemma is the Poincar\'{e}-Sobolev inequality of higher order in the elliptic setting.
\begin{lemma}\label{lemma: ce13}(\cite[Chapter 3]{L00}) \;\;Let $1\leq p<\infty,$ $R>0$ and $u$ be a smooth function in a ball $B:=B_{R}(x).$ Then there exists a polynomial $Q_{B} (u)$ of degree at most $m-1$ such that, for each $|\gamma|=0,1, 2..., m-1,$\begin{equation}\label{eq: cd1045} \fint_{B}D^{\beta}(u-Q_{B} (u))dx=0\end{equation}  and 
\begin{equation}\label{eq: ce14}\|D^{\gamma} (u-Q_{B} (u))\|_{L^{p}(B)}\leq C R^{m-|\gamma|} \|D^{m} u\|_{L^{p}(B)},\end{equation} where the constant $C=C(p)$ is independent of $R$ and $u.$
\end{lemma}

To build the H\"{o}lder estimates for solutions of parabolic equations of higher order, we need the lemma below. Its proof is quite similar to that in \cite[Lemma 4.3]{L1}, hence the details are omitted. 

\begin{lemma}\label{lemma: ce15} (\cite[Lemma 4.3]{L1})\;\; Assume $u\in L^{2}(Q_{2R}(t_{0}, x_{0}))$ and that there are positive constants $0<\eta \leq 1$ and $M$ along with a function $g$ defined on $Q_{2R}\times (0, R)$ such that $$\int_{Q(Y, r)}|u(X)-g(Y, r)|^{2}dX\leq M r^{n+2m+2\eta} $$ for any $Y=(s, y)\in Q_{R}$ and $r\in (0, R).$ Then, there exists a constant $C:=C(\eta)$ such that 
\begin{equation}\label{eq: ce16}[u]_{\eta, Q_{R}}\leq CM.\end{equation} 
\end{lemma}

For any $0\leq \lambda \leq n,$ we define the Morrey space $M^{2, \lambda} (\Omega)$ as the space of functions $f\in L^{2} (\Omega)$ such that $$\|f\|_{M^{2, \lambda} (\Omega)}:=\l(\sup_{x\in \Omega, \;\rho\leq 1}\rho^{-\lambda} \int_{B_{\rho}(x)}|f|^{2}\r)^{1/2}.$$ We stress that the two universal constants $C, c$ in the following contents inevitably depend on $\mu, l, \nu, c_{0}$ defined in Definition \ref{definition: ce17}.

\begin{lemma}\label{lemma: ce34} Assume that $\L\in \Ez_{2m}(c_{1}, 0, c_{2})$ has the H\"{o}lder property $H(\mu, \infty)$ as in Definition \ref{definition: ce17}. Suppose also that $u$ is a weak solution of the inhomogeneous elliptic equation  
\begin{equation}\label{eq: ce35} \L u=f \quad \mbox{in} \quad B_{2}:=B_{2}(x_{0}),\end{equation} where $f\in M^{2, \lambda}(B_{2})$ with $0\leq \lambda<n.$ Then for all $\ez\in (0, n-2m+2\mu)$ there exists a constant $C=C(\lambda,\ez)>0$ such that setting $$\varkappa:=\inf\{\lambda+2m, n-2m+2\mu-\ez\}<n$$ we have 
\begin{equation}\label{eq: ce36}\int_{B_{r}(x)}|\nabla^{m} u|^{2}\leq C r^{\varkappa}\l(\int_{B_{2}}|\nabla^{m} u|^{2}+\|f\|^{2}_{M^{2, \lambda}(B_{2})}\r)\end{equation} uniformly for all $x\in B_{1}(x_{0})$ and $0<r\leq 1.$ Hence, $u\in \Xi_{m-|\gamma|}^{2, \varkappa+2(m-|\gamma|)}(B_{1}),$ for all $0\leq |\gamma|\leq m,$ with the estimate 
\begin{equation}\label{eq: ce37}\|D^{\gamma}u\|_{\Xi_{m-|\gamma|}^{2, \varkappa+2(m-|\gamma|)}(B_{1})}\leq C \l(\|\nabla^{m}u\|_{L^{2}(B_{2})}+\|f\|_{M^{2, \lambda}(B_{2})}\r),\end{equation} where $$\|h\|_{\Xi_{m-|\gamma|}^{2, \varkappa+2(m-|\gamma|)}(B_{1})}:=\l(\sup_{x\in B_{1}(x_{0}), \;0\leq r\leq 1}\l[\rho^{-\varkappa-2(m-|\gamma|)}\inf_{P\in \P_{m-|\gamma|-1} }\int_{B_{\rho}(x)}|h(y)-P(y)|^{2}\r]\r)^{1/2}.$$ Here, for any integer $k\geq 1,$ $\P_{k-1}$ is the class of polynomials of degrees less than $k.$ \end{lemma}

{\it Proof.}\quad Note that, by \cite[Proposition 15]{AQ}, the H\"{o}lder property $H(\mu, \infty)$ of $\L$ implies that, for all $0<\rho<r<\infty$ and $x\in \rz,$  $$\int_{B_{\rho}(x)}|\nabla^{m} u|^{2}\leq c \l(\frac{\rho}{r}\r)^{n-2m+2\mu} \int_{B_{r}(x)}|\nabla^{m} u|^{2}.$$ In fact, this is the global Dirichlet property $D(\mu, \infty)$ in \cite[Definition 5]{AQ}. Moreover, \eqref{eq: ce36} is already established in \cite[Theorem 22]{AQ}. With \eqref{eq: ce36} in hand, by employing the Poincar\'{e}-Sobolev inequality in Lemma \ref{lemma: ce13} we can find a  polynomial $Q_{B}(u)$ of degree $m-1$ such that  
\begin{equation*} \int_{ B_{r} (x)}|D^{\gamma}u-D^{\gamma}Q_{B}(u)|^{2} \leq c r^{2(m-|\gamma|)}\int_{B_{r}(x)}|\nabla^{m} u|^{2}\leq c r^{2(m-|\gamma|)+\varkappa}\l(\int_{B_{2}}|\nabla^{m} u|^{2}+\|f\|^{2}_{M^{2, \lambda}(B_{2})}\r)
\end{equation*}
uniformly for all $x\in B_{1}(x_{0})$ and $0<r\leq 1,$ which immediately gives us the estimate \eqref{eq: ce37}.

\hfill$\Box$

\section{The H\"{o}lder property from elliptic equations to prabolic equations}

With the preparatory work from Section 3 in place, we are ready to prove our main results. Our first theorem is formulated as follows.

\begin{theorem}\label{theorem: ce20} If $\L\in \Ez_{2m}(c_{1}, 0, c_{2})$ satisfies the H\"{o}lder property $H(\mu, \infty)$ given in Defnition \ref{definition: ce17}, then the parabolic operator $\H$ with $t$-independent coefficients also satisfies the H\"{o}lder property $H(\mu_{0}, \infty)$ for some $\mu_{0}<\mu.$ The same statements hold for $\L^{*}$ and $\H^{*}.$ \end{theorem}


{\it Proof.}\quad Given any $(t_{0}, x_{0})\in \rr\times \rz$ and $R>0.$ Let $u$ be a weak solution of $\H$ in $Q_{R}:=[t_{0}-R^{2m}, t_{0}]\times B_{R}(x_{0}).$ Define $\tilde{u}(t, x):=u(t_{0}+R^{2m}t, x_{0}+Rx)$ and $\tilde{a}_{\alpha, \beta}(x)=\tilde{a}_{\alpha, \beta}(x_{0}+Rx).$ It is not hard to show that  $\tilde{u}(t, x)$ is a solution of $\tilde{\H}$ in $Q_{4}(0, 0),$ where 
$$\tilde{\H}:=\partial_{t}+\tilde{\L}:=\partial_{t}+\sum_{|\alpha|=|\beta|=m}(-1)^{|\alpha|}\partial^{\alpha}(\tilde{a}_{\alpha, \beta}(x)\partial^{\beta}).$$ Concerning the new elliptic operator $\tilde{\L}$, it has the H\"{o}lder property $H(\mu, \infty)$ if $\L$ has the H\"{o}lder property $H(\mu, \infty),$ see \cite{AQ}. Moreover, the elliptic constants $c_{1}, c_{2}$ remain the same for $\tilde{a}_{\alpha, \beta}(x).$ Thus, in the sequel, we can assume $R=4, t_{0}=0, x_{0}=0,$ otherwise we consider the parabolic operator $\tilde{\H}.$ 

Write $$\L u=-u_{t} \quad \mbox{in}\; Q_{4}:=Q_{4}(0, 0).$$ From \eqref{eq: ce28}, it follows that, for all $-2^{2m}\leq s \leq 0,$ \begin{equation}\label{eq: ce80}\|u_{t}(s, \cdot)\|_{L^{2}(B_{2})}\leq c \| u\|_{L^{2}(Q_{4})}.\end{equation} Clearly, the inequality \eqref{eq: ce80} enables us to apply Lemma \ref{lemma: ce34} with $f=-u_{t}(s, \cdot)$ and $\lambda=0$ ($M^{2, 0}(B_{2})=L^{2}(B_{2})$) to find a constant $$0<\varkappa:=\inf\{2m, n-2m+2\mu-\ez\}<n$$ such that the estimate \eqref{eq: ce37} holds. 

If $n\leq 2m,$ since $\mu=l+\nu$ and $|\gamma|\leq l,$ then we have $\varkappa+2(m-|\gamma|)\geq n+2\nu_{0}$ with $0<\nu_{0}<\nu<1$ by choosing $\ez$ small (for example, $\ez<\inf\{n-2m+2\mu, \nu\}$) when necessary. Note that, in light of  \cite{GM} or \cite{AQ}, we have, for any $k\geq 2,$ 
\begin{equation}\label{eq: y0}
L^{2}(B_{1}(0))\cap \Xi_{k}^{2, \lambda}(B_{1}(0))\backsimeq L^{2}(B_{1}(0))\cap \Xi_{k-1}^{2, \lambda}(B_{1}(0)) \quad \mbox{if}\; 0\leq \lambda<n+2k,\end{equation} and \begin{equation}\label{eq: y1}\Xi_{m-|\gamma|}^{2, \varkappa+2(m-|\gamma|)}(B_{1})\subset \Xi_{m-|\gamma|}^{2, n+2\nu_{0}}(B_{1}).\end{equation} Consequently, combining \eqref{eq: ce37} and \eqref{eq: y1} we get 
$$\|D^{\gamma}u\|_{\Xi_{m-|\gamma|}^{2, n+2\nu_{0} }(B_{1})}\leq C \l(\|\nabla^{m}u\|_{L^{2}(B_{2})}+\|f\|_{M^{2, \lambda}(B_{2})}\r).$$
From this and \eqref{eq: y0}, also recalling the definition of $\Xi_{1}^{2, n+2\nu_{0} }(B_{1})$ we see, 
for all $x\in B_{1}(0),$ all $0\leq |\gamma|\leq l,$ all $-2^{2m}\leq s \leq 0$ and all $0<r\leq 1,$ 
\begin{equation}\label{eq: ce38}\begin{aligned}\int_{B_{r}(x)}&|D^{\gamma} u(s, \cdot)-\fint_{B_{r}(x)}D^{\gamma} u(s, z)dz|^{2}\\[4pt]
&\leq c r^{n+2\nu_{0} }\l(\int_{B_{2}}|D^{m} u(s, \cdot)|^{2}+\|u_{t}(s, \cdot)\|^{2}_{L^{2}(B_{2})}+\int_{B_{2}}|D^{\gamma} u(s, \cdot)|^{2}\r)\\[4pt]&\leq c r^{n+2\nu_{0} }\|u\|^{2}_{L^{2}(Q_{4})},
\end{aligned}\end{equation} where in the last step we used \eqref{eq: ce27}-\eqref{eq: ce28}. 

Then, by \eqref{eq: ce38}, we have that, for all $(t, x)\in Q_{1},$ $0<r\leq 1$ and $0\leq |\gamma|\leq l,$ \begin{equation}\label{eq: ce39}\begin{aligned}
\inf_{a\in \rr}\int_{t-r^{2m}}^{t}\int_{B_{r}(x)}|D^{\gamma}u-a |^{2}dyds &\leq r^{2m} \inf_{a\in \rr}\sup_{s\in (t-r^{2m}, t) }\int_{B_{r}(x)}|D^{\gamma}u(s, \cdot)-a |^{2}dyds\\[4pt]
& \leq c r^{2m+n+2\nu_{0}}\|u\|^{2}_{L^{2}(Q_{4})}.
\end{aligned}
\end{equation}  Thus, by \eqref{eq: ce39},  there exists a constant $L,$ depending on $u, r, x,$ such that 
\begin{equation}\label{eq: ce40}\begin{aligned}
\int_{t-r^{2m}}^{t}\int_{B_{r}(x)}|D^{\gamma}u(s, y)-L|^{2}dsdy &\leq c r^{2m+n+2\nu_{0}}\|u\|^{2}_{L^{2}(Q_{4})}
\end{aligned}
\end{equation} for all $(t, x)\in Q_{1}$ and $0<r\leq 1,$ where $0<\nu_{0}<\nu<1.$ Now, by \eqref{eq: ce40}, an application of Lemma \ref{lemma: ce15} leads to 
\begin{equation}\label{eq: ce41}
[D^{\gamma} u]_{\nu_{0}, Q_{1}(0, 0)}\leq c\|u\|_{L^{2}(Q_{4})},\quad \forall \; 0\leq |\gamma|\leq l.
\end{equation} 
Moreover, using \eqref{eq: ce41} and the argument at the beginning (translations and dilatations), we can conclude that for any $(t_{0}, x_{0})\in \rr\times \rz$ and $R>0,$ if $u$ is a weak solution of $\H$ in $Q_{4R}:=[t_{0}-(4R)^{2m}, t_{0}]\times B_{R}(x_{0}),$ then $$[D^{\gamma} u]_{\nu_{0}, Q_{R}(t_{0}, x_{0})}\leq c R^{-(\frac{n}{2}+m+|\gamma|+\nu_{0})}\|u\|_{L^{2}(Q_{4R}(t_{0}, x_{0}))},\quad \forall \;\; 0\leq |\gamma|\leq l.$$

In the case $2m+1\leq n\leq 2m+2l+1,$ we argue as follows. For clarity, we now suppose that $u$ is a solution of $\H$ in $Q_{8}(0, 0).$ Fix $(t_{0}, x_{0}) \in Q_{2}(0, 0).$ Clearly, $u_{t}$ is a solution of $\H$ in $Q_{4}(t_{0}, x_{0})$ thanks to the $t$-independence of the coefficients. For any $s\in [t_{0}-2^{2m}, t_{0}],$ using Lemma \ref{lemma: ce26} and \eqref{eq: ce23}, we see $$\|u_{tt}(s, \cdot)\|^{2}_{L^{2}(B_{2}(x_{0}))}\lesssim \|u_{t}\|^{2}_{L^{2}(Q_{4}(t_{0}, x_{0}))}\lesssim \|u\|^{2}_{L^{2}(Q_{8}(0, 0))}.$$ Since $\L u_{t} =u_{tt},$ then, applying Lemma \ref{lemma: ce34}, also Lemma \ref{lemma: ce26} and \eqref{eq: ce23}, we obtain a constant $$\varkappa:=\inf\{2m, n-2m+2\mu-\ez\}<n$$ such that for all $x\in B_{1}(x_{0})$ and all $0<r\leq 1,$ 
\begin{equation}\label{eq: ce42}\begin{aligned}
\|u_{t}(s, \cdot)\|^{2}_{\Xi_{m}^{2, \varkappa+2m}(B_{1}(x_{0}))}&\lesssim \l(\int_{B_{2}(x_{0})}|\nabla^{m} u_{t}(s, \cdot)|^{2}+\|u_{tt}(s, \cdot)\|^{2}_{L^{2}(B_{2}(x_{0}))}\r)\\[4pt]
&\lesssim \|u\|^{2}_{L^{2}(Q_{8}(0, 0))}\end{aligned}\end{equation} uniformly in $s\in [t_{0}-2^{2m}, t_{0}].$
Since \eqref{eq: ce42} holds for all $(t_{0}, x_{0}) \in Q_{2}(0)$ and $s\in [t_{0}-2^{2m}, t_{0}],$ we get that 
$$\|u_{t}(s, \cdot)\|^{2}_{\Xi_{m}^{2, \varkappa+2m}(B_{2}(0))} \leq c \|u\|^{2}_{L^{2}(Q_{8} (0, 0) )} \quad \mbox{uniformly in} \; s\in [-2^{2m}, 0].$$
Again, according to \cite{GM} (see also \cite{AQ}), we have that $$\Xi_{m}^{2, \varkappa+2m}(B_{2}(0) )\subset \Xi_{m}^{2, \varkappa}(B_{2}(0))$$ and $$L^{2}(B_{2}(0))\cap \Xi_{m}^{2, \varkappa}(B_{2}(0))\backsimeq L^{2}(B_{2}(0))\cap \Xi_{1}^{2, \varkappa}(B_{2}(0) )\subset M^{2, \varkappa}(B_{2}(0))$$ since $\varkappa<n.$ Thus, for any $s\in [-2^{2m}, 0],$ exploiting \eqref{eq: ce28} and the above norm equivalence, we arrive at
\begin{equation}\label{eq: ce44}\begin{aligned}
\|u_{t}(s, \cdot)\|_{M^{2, \varkappa}(B_{2}(0))}&\leq c \|u_{t}(s, \cdot)\|_{L^{2}(B_{2}(0))}+c \|u_{t}(s, \cdot)\|_{\Xi_{m}^{2, \varkappa}(B_{2}(0))} \\[4pt]
&\leq c\|u\|^{2}_{L^{2}(Q_{8} (0, 0) )}+c \|u_{t}(s, \cdot)\|_{\Xi_{m}^{2, \varkappa+2m}(B_{2}(0))}\\[4pt]
&\leq c \|u\|^{2}_{L^{2}(Q_{8} (0, 0) )}.
\end{aligned}
\end{equation}



As a result, we are allowed by the estimate \eqref{eq: ce44} to apply Lemma \ref{lemma: ce34} with $f=-u_{t}(s, \cdot)$ and $\lambda=\varkappa_{0}$-we can choose $\varkappa_{0}=1+2l$ since $\varkappa\geq 1+2l$ and $l\leq m-1$-to obtain a positive constant $$\varkappa_{1}:=\inf\{2m+\varkappa_{0}, n-2m+2\mu-\ez\}=n-2m+2\mu-\ez<n \;(n\leq 2m+2l+1)$$ such that \begin{equation}\label{eq: y5}\|D^{\gamma}u\|_{\Xi_{m-|\gamma|}^{2, \varkappa_{1}+2(m-|\gamma|)}(B_{1}(0))}\leq C \l(\|\nabla^{m}u\|_{L^{2}(B_{2}(0))}+\|f\|_{M^{2, \lambda}(B_{2}(0))}\r).\end{equation} Note that $\varkappa_{1}+2(m-|\gamma|)=n+2(l-|\gamma|)+2\nu-\ez\geq n+2\nu_{0},$ with $0<\nu_{0}<\nu<1,$ by taking $\ez$ small enough. A similar argument as in \eqref{eq: y0}-\eqref{eq: ce41} gives us 
\begin{equation}\label{eq: ce47}
[D^{\gamma} u]_{\nu_{0}, Q_{1}(0, 0)}\leq c \|u\|_{L^{2}(Q_{8}(0, 0))},\quad \forall \; 0\leq |\gamma|\leq l.
\end{equation}


If $n\geq 2(m+l+1),$ we repeat the above procedure. Apparently, the process cannot go on indefinitely, so, by \eqref{eq: ce47}, $$ [D^{\gamma} u]_{\nu_{0}, Q_{1}(0, 0)}\leq c\|u\|_{L^{2}(Q_{8}(0, 0))},\quad \forall \; 0\leq |\gamma|\leq l,$$ is concluded. Armed with this estimate, by translations and dilatations, also a well known covering argument, we are easy to obtain, for some $ 0<\nu_{0}<\nu<1,$ \begin{equation}\label{eq: ce48}[D^{\gamma} u]_{\nu_{0}, Q_{r}(t_{0}, x_{0})}\leq cr^{-(\frac{n}{2}+m+|\gamma|+\nu_{0})}\|u\|_{L^{2}(Q_{2r}(t_{0}, x_{0}))},\; \forall \; \;0\leq |\gamma|\leq l.\end{equation}

Fix $(s, y)\in Q_{r}(t_{0}, x_{0}).$ Then, for all $(\tau, \xi)\in Q_{r}(s, y) \subset Q_{2r}(t_{0}, x_{0}),$ it follows from \eqref{eq: ce48} that \begin{equation*}\begin{aligned}
|D^{\gamma} u(s, y)|&\leq |D^{\gamma} u(s, y)-D^{\gamma} u(\tau, \xi)|+|D^{\gamma} u(\tau, \xi)|\\[4pt]
&\leq (|\tau-s|^{\frac{1}{2m}}+|y-\xi|)^{\nu_{0}}[D^{\gamma} u]_{\nu_{0}, Q_{r}(t_{0}, x_{0})}+|D^{\gamma} u(\tau, \xi)|\\[4pt]
&\leq c r^{-(\frac{n}{2}+m+|\gamma|)}\|u\|_{L^{2}(Q_{2r}(t_{0}, x_{0}))}+|D^{\gamma} u(\tau, \xi)|.\end{aligned}
\end{equation*}
 By averaging the above inequality over $ Q_{r}(s, y)$ with respect to $(\tau, \xi),$ we see 
\begin{equation*}\begin{aligned}
|D^{\gamma} u(s, y)|&\leq c r^{-(\frac{n}{2}+m+|\gamma|)}\|u\|_{L^{2}(Q_{2r}(t_{0}, x_{0}))}+c\fint_{Q_{r}(s, y)}|D^{\gamma} u(\tau, \xi)|d\tau d\xi\\[4pt]
&\leq c r^{-(\frac{n}{2}+m+|\gamma|)}\|u\|_{L^{2}(Q_{2r}(t_{0}, x_{0}))}+c\l(\fint_{Q_{r}(s, y)}|D^{\gamma} u(\tau, \xi)|^{2}d\tau d\xi\r)^{1/2} \\[4pt]
&\leq c r^{-(\frac{n}{2}+m+|\gamma|)}\|u\|_{L^{2}(Q_{2r}(t_{0}, x_{0}))} +cr^{-|\gamma|}\l(\fint_{Q_{2 r}(s, y)}|u(\tau, \xi)|^{2}d\tau d\xi\r)^{1/2} \quad (\mbox{by \eqref{eq: ce5}}) \\[4pt]
&\leq c r^{-(\frac{n}{2}+m+|\gamma|)}\|u\|_{L^{2}(Q_{3r}(t_{0}, x_{0}))},
\end{aligned}
\end{equation*} for any $0\leq |\gamma|\leq l.$ Therefore, \begin{equation}\label{eq: ce50}r^{|\gamma|}\sup_{(s, y)\in Q_{r}(t_{0}, x_{0})}|D^{\gamma} u(s, y)|\leq c \l(\fint_{Q_{2 r}(t_{0}, x_{0}))}|u(\tau, \xi)|^{2}d\tau d\xi\r)^{1/2},\quad \forall \; 0\leq |\gamma|\leq l,\end{equation} holds by a covering argument again.  If we set $\mu_{0}:=l+\nu_{0},$ then $\mu_{0}<\mu,$ and \eqref{eq: ce48} with $|\gamma|=l$ and \eqref{eq: ce50} yield that $\H$ has the H\"{o}lder property $H(\mu_{0}, \infty).$ This suffices.


\hfill$\Box$

\begin{remark}\label{remark: ce49} If $u$ is a solution of $\H$ in $Q_{R}(t, x),$ we can conclude from \eqref{eq: ce50} that for any $Q_{2 r}(t_{0}, x_{0}) \subset Q_{R}(t, x),$
$$\sup_{(s, y)\in Q_{3/2r}(t_{0}, x_{0})}| u(s, y)|\leq c \l(\fint_{Q_{2 r}(t_{0}, x_{0}))}|u(\tau, \xi)|^{2}d\tau d\xi\r)^{1/2},$$ which implies, for any $q>2,$ \begin{equation}\label{eq: ce51}\l(\fint_{Q_{ 3/2r}(t_{0}, x_{0}))}|u(\tau, \xi)|^{q}d\tau d\xi\r)^{1/q}\leq c \l(\fint_{Q_{2 r}(t_{0}, x_{0}))}|u(\tau, \xi)|^{2}d\tau d\xi\r)^{1/2}.\end{equation} From \eqref{eq: ce51} and the self-improved lemma \cite[Theorem B.1]{BCF}, it follows that \begin{equation}\label{eq: ce52}\l(\fint_{Q_{ 3/2r}(t_{0}, x_{0}))}|u(\tau, \xi)|^{q}d\tau d\xi\r)^{1/q}\leq c \fint_{Q_{2 r}(t_{0}, x_{0}))}|u(\tau, \xi)|d\tau d\xi.\end{equation} The estimate \eqref{eq: ce52} contributes to that \eqref{eq: ce50} and \eqref{eq: ce48} with $|k|=l$ can be improved to $$r^{|\gamma|}\sup_{(s, y)\in Q_{r}(t_{0}, x_{0})}|D^{\gamma} u(s, y)|\leq c \fint_{Q_{2 r}(t_{0}, x_{0}))}|u(\tau, \xi)|d\tau d\xi,\quad \forall \; 0\leq |\gamma|\leq l,$$ and for $|\gamma|=l,$ $$[D^{\gamma} u]_{\nu_{0}, Q_{r}(t_{0}, x_{0})}\leq c r^{-(l+\nu_{0})}\fint_{Q_{2 r}(t_{0}, x_{0}))}|u(\tau, \xi)|d\tau d\xi,$$ respectively.

\end{remark}

As a straightforward consequence of Theorem \ref{theorem: ce20}, we have the following perturbation theorem. A result of this type mentioned in the introduction for second order parabolic equations was essentially established in \cite{A}. 
\begin{theorem}\label{theorem: ce59}\; Assume that $\H_{0}:= \partial_{t} +\L_{0} :=\partial_{t} +\sum_{|\alpha|=|\beta|=m}(-1)^{|\alpha|}\partial^{\alpha}(a^{0}_{\alpha, \beta}(x)\partial^{\beta})$ and $\H_{1}:=\partial_{t} +\L_{1} :=\partial_{t} +\sum_{|\alpha|=|\beta|=m}(-1)^{|\alpha|}\partial^{\alpha}(a^{1}_{\alpha, \beta}(x)\partial^{\beta})$ satisfy \eqref{eq: ce200} and \eqref{eq: ce1}, and in addition, $\H_{0}$ has the H\"{o}lder property $H(\mu, \infty).$ If $$\sup_{|\alpha|=|\beta|=m}\|a^{0}_{\alpha, \beta}-a^{1}_{\alpha, \beta}\|_{\infty} <\ez,$$ then there exist a constant $\ez_{0}$ ($0<\ez<\ez_{0}$) and some $\mu_{0}<\mu$ such that $\H_{1}$ verifies the H\"{o}lder property $H(\mu_{0}, \infty).$ 

\end{theorem}
{\it Proof.}\quad Since the coefficients are $t-$independent, then $\L_{0}$ satisfies the H\"{o}lder property $H(\mu, \infty)$ (\eqref{eq: ce18}-\eqref{eq: ce19}). Hence, invoking \cite[Lemma 42]{AQ} we immediately get that $\L_{1}$ has the H\"{o}lder property $H(\mu_{1}, \infty)$ for any $0<\mu_{1}<\mu$ when $\ez$ small enough. This allows us to employ Theorem \ref{theorem: ce20} to find a constant $\mu_{0}<\mu_{1}$ such that $\H_{1}$ verifies $H(\mu_{0}, \infty),$ that is, \eqref{eq: ce48} with $|\gamma|=l$ and \eqref{eq: ce50} associated to $\mu_{0}$. The proof is finished. 

\hfill$\Box$

\begin{remark}\label{remark: ce60}  Based on Theorem \ref{theorem: ce59}, if $\{a^{0}_{\alpha, \beta}\}$ are either complex constants or functions belonging to the space $VMO$-We say that a function $f\in BMO$ belongs to $VMO$ if $$\lim_{\rho\to 0}\sup_{B_{r}, 0<r\leq \rho} \l(\fint_{B_{r}}|f-\fint_{B_{r}}f|^{2}\r)^{1/2}=0.$$-and $\{a^{1}_{\alpha, \beta}\}$ are as in Theorem \ref{theorem: ce59}, then both $\H_{1}$ and $\H^{*}_{1}$ verify the H\"{o}lder property $H(\mu_{0}, \infty)$ by a combination of \cite[Proposition 45]{AQ}, \cite[Proposition 50]{AQ} and \cite[Theorem 12]{AQ}.

\end{remark}

\section{Gaussian estimates for the fundamental solution of $\H$ }

Gaussian upper bounds and H\"{o}lder regularity estimates (\eqref{eq: ce3002}-\eqref{eq: ce3004}) of the heat kernel of the semigroup $e^{-t\L}$ can be derived from the H\"{o}lder properties $H(\mu, \infty)$ of $\L$ and its adjoint in light of \cite[Theorem 12]{AQ}. This fact, coupled with the $t-$independence of the coefficients, implies that the fundamental solution of $\H$ constructed by using the same way as that in \cite{NAS, NK} satisfies similar Gaussian estimates under the $H(\mu, \infty)$-assumptions for $\H$ and $\H^{*}$. Here, we derive pointwise bounds for the fundamental solution of $\H,$ which differs from the statements in \cite[Theorem 3.4, Corollary 3.6]{K}. These bounds primarily serve as essential preparation for extending the results of \cite{NAS, NK}  to higher-order cases.


\begin{theorem}\label{theorem: ce52}\; Suppose that both $\H$ and $\H^{*}$ with $\L \in \Ez_{2m}(c_{1}, 0, c_{2})$ satisfy the H\"{o}lder property $H(\mu, \infty).$ Then $\H$ has a fundamental solution $\Gamma(t, x, s, y)$ such that for all $x, y, h\in \rz,$ all $t>s$ and all multi-index $\gamma\in \nn^{n},$
\begin{equation}\label{eq: ce53}
|D_{x}^{\gamma}\Gamma(t, x, s, y)|+|D_{y}^{\gamma}\Gamma(t, x, s, y)|\leq \frac{C}{(t-s)^{\frac{n+|\gamma|}{2m}}}g_{m, c}\l(\frac{|x-y|}{(t-s)^{\frac{1}{2m}}}\r),\end{equation} when $|\gamma|\leq l,$
\begin{equation}\label{eq: ce54}
|D_{x}^{\gamma}\Gamma(t, x+h, s, y)-D_{x}^{\gamma}\Gamma(t, x, s, y)|\leq \frac{C}{(t-s)^{\frac{n+|\gamma|}{2m}}}\l(\frac{|h|}{(t-s)^{\frac{1}{2m}}+|x-y|}\r)^{\nu}g_{m, c}\l(\frac{|x-y|}{(t-s)^{\frac{1}{2m}}}\r)\end{equation} 
and \begin{equation}\label{eq: ce55}
|D_{y}^{\gamma}\Gamma(t, x, s, y+h)-D_{y}^{\gamma}\Gamma(t, x, s, y)|\leq \frac{C}{(t-s)^{\frac{n+|\gamma|}{2m}}}\l(\frac{|h|}{(t-s)^{\frac{1}{2m}}+|x-y|}\r)^{\nu}g_{m, c}\l(\frac{|x-y|}{(t-s)^{\frac{1}{2m}}}\r),
\end{equation} when $|\gamma|=l$ and $2|h|\leq (t-s)^{\frac{1}{2m}}+|x-y|.$ In particular, \begin{equation}\label{eq: ce56}
\Gamma(t, x, s, y)=0\;\mbox{if}\; t\leq s\; \;\mbox{and}\; \;\Gamma(t, x, s, y)=\Gamma(t-s, x, 0, y).\end{equation} Let $\Gamma^{*}(t, x, s, y)$ denote the fundamental solution of $\H^{*},$ then 
\begin{equation}\label{eq: ce57}
\Gamma^{*}(t, x, s, y)=\overline{\Gamma(s, y, t, x)},\;\quad  \forall \; t>s. 
\end{equation} Obviously, $\Gamma^{*}$ has the similar properties as $\Gamma.$

\end{theorem}
{\it Proof.}\quad Since the coefficients are $t-$independent and $\H$ has the H\"{o}lder property $H(\mu, \infty)$, it follows that any solution $u$ of $\L$ satisfies \eqref{eq: ce18}-\eqref{eq: ce19} associated to $\mu.$ The same conclusion holds for $\L^{*}.$ Hence, by \cite[Theorem 12]{AQ}, the heat kernel $K_{t}(x, y)$ of $\L$ satisfies the Gaussian estimates defined in Definition \ref{definition: ce3001} for any $0<t<\infty$. If we define $$\Gamma(t, x, s, y):=K_{t-s}(x, y), \;\mbox{if}\; t>s,\; \Gamma(t, x, s, y):=0, \;\mbox{if}\; t\leq s,$$
then, $\Gamma$ is a fundamental solution of $\H$ in $\rdm.$ Moreover, \eqref{eq: ce53}-\eqref{eq: ce57} follows easily from the properties of $K_{t}(x, y)$ and its adjoint, see \cite{AQ}.

\hfill$\Box$
\begin{remark}\label{remark: ce101} \eqref{eq: ce53}-\eqref{eq: ce55} are indeed equivalent to that for all $x, y, h\in \rz$ and $t>s,$
\begin{equation}\label{eq: ce6001}|D_{x}^{\gamma}\Gamma(t, x, s, y)|+|D_{y}^{\gamma}\Gamma(t, x, s, y)|\leq \frac{C}{(t-s)^{\frac{n+|\gamma|}{2m}}}g_{m, c}\l(\frac{|x-y|}{(t-s)^{\frac{1}{2m}}}\r),
\end{equation}when $|\gamma|\leq l,$\begin{equation}\label{eq: ce6002}
|D_{x}^{\gamma}\Gamma(t, x+h, s, y)-D_{x}^{\gamma}\Gamma(t, x, s, y)|\leq \frac{C}{(t-s)^{\frac{n+|\gamma|}{2m}}}\l(\frac{|h|}{(t-s)^{\frac{1}{2m}}}\r)^{\nu}\end{equation}
and \begin{equation}\label{eq: ce6003}
|D_{y}^{\gamma}\Gamma(t, x, s, y+h)-D_{y}^{\gamma}\Gamma(t, x, s, y)|\leq \frac{C}{(t-s)^{\frac{n+|\gamma|}{2m}}}\l(\frac{|h|}{(t-s)^{\frac{1}{2m}}}\r)^{\nu},\end{equation} when $|\gamma|= l.$ It is more convenient to prove \eqref{eq: ce6001}-\eqref{eq: ce6003} but \eqref{eq: ce53}-\eqref{eq: ce55} are more useful in pratice. 
\end{remark}

\begin{remark}\label{remark: ce100} A local version of Theorem \ref{theorem: ce52} (corresponding to the local Definition \ref{definition: ce3001}) can be obtained in the same way by introducing the local H\"{o}lder property $H(\mu, R_{0})$ in \cite[Definition 6]{AQ}. Of course, the constants $C, c, c_{0}$ in Definition \ref{definition: ce17} and and Theorem \ref{theorem: ce52} will depend on $R_{0}.$ 
\end{remark}

\begin{remark}\label{remark: ce58}  In veiw of the theory of second order parabolic equations (systems) \cite{ AN, AE, HK}, Theorem \ref{theorem: ce52} shall hold for parabolic equations (systems) of higher order with time-dependent coefficients, even in a weighted setting. This will be the focus of our future work.
\end{remark}

In particular, when $2m\geq n,$ the assumption that $\H$ and $\H^{*}$ have the H\"{o}lder property $H(\mu, \infty)$ in Theorem \ref{theorem: ce52} can be removed.  

\begin{theorem}\label{theorem: ce1010}\; The fundamental solution $\Gamma(t, x, s, y)$ of any parabolic operator $\H$ defined in \eqref{eq: ce00} with $\L\in \Ez_{2m}(c_{1}, 0, c_{2})$ satisfies \eqref{eq: ce53}-\eqref{eq: ce55} with an integer $l$ and some $0<\nu<1$.
\end{theorem}
{\it Proof.}\quad We begin with two simple facts. First, by a similar construction in \cite[Remark, Section 0.2]{AT}, $H^{-m}(\rz)$ is the set of vectors of the form $u=\sum_{|\alpha|\leq m}\partial^{\alpha} u_{\alpha}$ with $u_{\alpha} \in L^{2}(\rz),$ and its norm is defined by $\|u\|_{H^{-m}}:=\inf\{\l(\sum_{|\alpha|\leq m}\|u_{\alpha}\|_{2}^{2}\r)^{1/2}\}$ where the infimum is taken over all representations of $u.$ Second, for any $\L \in \Ez_{2m}(c_{1}, \lambda_{0}, c_{2}, c_{3}),$ it follows from \cite[Lemma 3]{AQ} that there exists a $\lambda:=c_{1}/2+\lambda_{0}+c_{3}$ such that $\L +\lambda$ is maximal accretive of type $\omega<\frac{\pi}{2}$ and is the generator of a contraction semigroup $e^{-t(\L+\lambda)}$ that is analytic. Moreover, imitating the proof of \cite[Proposition 1 Chapter 1]{AT}, we can conclude that, for a fixed $\theta\in (\omega, \frac{\pi}{2})$ and all $z\in \Gamma_{\pi-\theta},$ all $u\in H^{m}$ and $|\gamma|\leq m,$ $$|z|^{\frac{|\gamma|}{2m}}\|D^{\gamma}(z+\L+\lambda)^{-1}u\|_{2}+|z|^{\frac{|\gamma|}{2m}}\|(z+\L+\lambda)^{-1}D^{\gamma}u\|_{2}\leq C \|u\|_{2}$$ and $$|z|^{\frac{m+|\gamma|}{2m}}\|D^{m}(z+\L+\lambda)^{-1}D^{\gamma}u\|_{2}\leq C \|u\|_{2}$$ where $C$ depends only on $c_{1}, c_{2}, c_{3}, \omega, \theta, \lambda_{0}.$ The two facts together imply that, for any $z\in \Gamma_{\pi-\theta}$ with $|z|=1,$ \begin{equation}\label{eq: ce201}(z+\L+\lambda)^{-1}: H^{-m}\to H^{m}\quad \mbox{ is bounded}\end{equation} with the operator norm depends only on $c_{1}, c_{2}, c_{3}, \omega, \theta, \lambda_{0}.$ 

\noindent{\textbf{Case 1: 2m>n.}} If $2m>n,$ the Sobolev embedding theorem gives us $$H^{m}(\rz)\subset C^{l, \nu}(\rz)$$ where $l=m-[\frac{n}{2}]-1$ and $0<\nu<1$. By this and the imbedding $H^{m}\subset L^{2},$ \eqref{eq: ce201} contributes to, for any $z\in \Gamma_{\pi-\theta}$ with $|z|=1,$ \begin{equation}\label{eq: ce202}\|(z+\L+\lambda)^{-1} f\|_{C^{l, \nu}}\leq C \|f\|_{2}.\end{equation} In particular, we can choose $z=c_{1}/2$ in the above estimates such that \begin{equation}\label{eq: ce203}\|(\L+c_{1}+\lambda_{0}+c_{3})^{-1} f\|_{C^{l, \nu}}\leq C \|f\|_{2}\end{equation} holds for all $\L \in \Ez_{2m}(c_{1}, \lambda_{0}, c_{2}, c_{3}).$ 

In the sequel, we consider $\L\in \Ez_{2m}(c_{1}, 0, c_{2}),$ that is $\lambda_{0}=0$ and $c_{3}=0.$ Thus, for any $z\in \Gamma_{\pi-\theta}$ with $|z|=1,$ \begin{equation}\label{eq: ce204}\|(z+\L)^{-1} f\|_{C^{l, \nu}}\leq C \|f\|_{2}\end{equation}  with $C$ being independent of $z.$ Set $Tf:=T_{s}f:=f(x/s)$ for $s>0.$ An calculation yields $$T^{-1}\L T=s^{-2m}\L^{s}$$ where $\L^{s}=\sum_{|\alpha|=|\beta|=m}(-1)^{|\alpha|}\partial^{\alpha}(a_{\alpha, \beta}(sx)\partial^{\beta}).$ Clearly, $\L^{s}\in \Ez_{2m}(c_{1}, 0, c_{2})$ and $T^{-1} (\L+z)^{-1}T=s^{2m}(s^{2m}z+\L^{s}).$ Applying \eqref{eq: ce204} with $\L$ replaced by $\L^{s}$ we get, for any $|\gamma|\leq l$ and all $z\in \Gamma_{\nu}$ with some $\theta \in (\frac{\pi}{2}, \pi-\omega),$\begin{equation}\label{eq: ce205}\|D^{\gamma}(z+\L)^{-1} f\|_{L^{\infty}}\leq C|z|^{\frac{n}{4m}+\frac{|\gamma|}{2m}-1} \|f\|_{2}\end{equation} and \begin{equation}\label{eq: ce206}[[(z+\L)^{-1} f]]_{\nu, \rz}\leq C|z|^{\frac{n}{4m}+\frac{l+\nu}{2m}-1} \|f\|_{2}\end{equation} by taking $s=|z|^{-\frac{1}{2m}}.$ With \eqref{eq: ce205}-\eqref{eq: ce206} in hand, using the Cauchy formular $$e^{-t\L}u=\frac{1}{2\pi i}\int_{\Sigma_{\pm}\cup \Sigma_{0}}e^{tz}(z+\L)^{-1}udz,$$ where $\Sigma_{\pm}:=\{z=re^{i\pm \theta}, r\geq \frac{1}{t}\}$ and $\Sigma_{0}:=\{z=\frac{1}{t}e^{i\eta}: |\eta|\leq \theta\},$ we can deduce \begin{equation}\label{eq: ce207}\|D^{\gamma}e^{-t\L}u\|_{\infty}\leq C t^{-\frac{n}{4m}-\frac{|\gamma|}{2m}}  \|u\|_{2}\quad \forall\; |\gamma|\leq l
\end{equation} and \begin{equation}\label{eq: ce208}
[[D^{l}e^{-t\L}u]]_{\nu, \rz}\leq Ct^{-\frac{n}{4m}-\frac{l+\nu}{2m}} \|u\|_{2}.
\end{equation} 

Let us admit \eqref{eq: ce3002} for the moment.  Then, we only have to show that \begin{equation}\label{eq: ce209}\sup_{y\in \rz}[[D_{x}^{\gamma}K_{t}(\cdot, y)]]_{\nu, \rz}\leq C t^{-\frac{n}{2m}-\frac{l+\nu}{2m}}, \quad  \forall\; |\gamma|= l,
\end{equation} in order to build \eqref{eq: ce3003}. Furthermore, \eqref{eq: ce209} is equivalent to $$[[D^{l}e^{-t\L}u]]_{\nu, \rz}\leq Ct^{-\frac{n}{2m}-\frac{l+\nu}{2m}} \|u\|_{1}$$ by \cite[Lemma 17 Chapter 1]{AT}. Note that \eqref{eq: ce207} holds with $\L^{*}.$ Then, by duality, we see that \begin{equation}\label{eq: ce2020}\|e^{-t\L}u\|_{L^{2}}\leq Ct^{-\frac{n}{4m}}\|u\|_{L^{1}}.\end{equation} This inequality and \eqref{eq: ce208} lead to $$[[D^{l}e^{-t\L}u]]_{\nu, \rz}\leq Ct^{-\frac{n}{2m}-\frac{l+\nu}{2m}} \|e^{-t/2\L}u\|_{2}\leq Ct^{-\frac{n}{2m}-\frac{l+\nu}{2m}} \|u\|_{L^{1}}.$$ This proves \eqref{eq: ce3003}. Besides, repeating the above arguments for $\L^{*}$ we can also deduce \eqref{eq: ce3004}. 

Next we seek to prove \eqref{eq: ce3002}. To the end, we use an idea stemming from \cite{D0} which is also exploited in \cite{AQ, AT}. Let $\phi \in C_{0}^{\infty}(\rz)$ be real-valued function such that $$\|D^{\gamma}\phi\|_{\infty}\leq 1,\quad \mbox{for all}\; 1\leq |\gamma|\leq m.$$ Define the operator of multiplication $e^{\phi}$ and consider $K_{t}^{\ez\phi}:=e^{-\ez\phi}e^{-t\L}e^{\ez\phi}$ for $t>0$ and $\ez>0.$ If we can prove for $\ez\geq c$ \begin{equation}\label{eq: ce2011}\|D^{\gamma}K_{1}^{\ez\phi}u\|_{\infty}\leq C \P(\ez)  e^{a\ez^{2m}}\|u\|_{1}\quad  \forall\; |\gamma|\leq l,
\end{equation} where $C, a>0$ are independent of $\ez$ and $\phi,$ and $\P(\ez):=\ez^{|\gamma|+n},$ then \begin{equation}\label{eq: ce2012}|D^{\gamma}K_{1}(x, y)|\leq C e^{-b|x-y|^{\frac{2m}{2m-1}}},\quad  \forall\; |\gamma|\leq l.
\end{equation} Indeed, If $|x-y|\leq c,$ \eqref{eq: ce2012} follows from \eqref{eq: ce207}and \eqref{eq: ce2020}. If $|x-y|\geq c,$ \eqref{eq: ce2012} can be attributed to \eqref{eq: ce2011} for $\ez\geq c$ according to the choice of $\ez$ in the following argument. To explain this point, we first note that the kernel of $K_{t}^{\ez\phi}$ is $e^{-\ez\phi(x)} K_{t}(x, y)e^{-\ez\phi(y)}$ since the kernel of $e^{-t\L}$ is $K_{t}(x, y).$ Invoking \cite[Lemma 17 Chapter 1]{AT} again, we get $$|D_{x}^{\gamma}(e^{-\ez\phi(x)} K_{1}(x, y)e^{-\ez\phi(y)}|\leq C \P(\ez)  e^{a\ez^{2m}}, \;\forall \; |\gamma|\leq l.$$ By Leibniz's rule, a tedius computation leads to $$|D^{\gamma}K_{1}(x, y)|\leq C(1+\ez+...+\ez^{|\gamma|})\P(\ez)e^{a\ez^{2m}}e^{\ez(\phi(x)-\phi(y))}.$$ We further choose $\phi$ such that $\phi(x)-\phi(y)=\frac{|x-y|}{2}$ (see \cite{AQ}), then \eqref{eq: ce2012} is easy to be derived from optimizing with respect to $\ez\geq c$ (In fact, $\ez=\frac{|x-y|^{\frac{1}{2m-1}} }{c(a)}$). Then, by scalling above, we can apply \eqref{eq: ce2012} to attain $$ |D^{\gamma}K^{\L^{1/2m}}_{1}(x, y)|\leq C e^{-b|x-y|^{\frac{2m}{2m-1}}}\quad \forall\; |\gamma|\leq l.$$ From this estimate and the fact that $K_{t}(x, y)=t^{-\frac{n}{2m}}K_{1}^{\L^{1/2m}}(\frac{x}{t^{1/2m}}, \frac{x}{t^{1/2m}}),$ it follows that 
\begin{equation}\label{eq: ce2013}
|D^{\gamma}_{x}K_{t}(x, y)|\leq Ct^{-\frac{n+|\gamma|}{2m}} e^{-c\l(\frac{|x-y|^{2m}}{t} \r)^{\frac{1}{2m-1}}},
\end{equation} that is \eqref{eq: ce3002}. Hence, we are left to build \eqref{eq: ce2011}. 

Since $\L \in \Ez_{2m}(c_{1}, 0, c_{2}),$ we then see $L^{\phi}:=e^{-\phi}e^{-\L}e^{\phi} \in \Ez_{2m}(c_{1}, 0, c_{2}, c_{3}).$ Indeed, $c_{3}$ only depends on $m$ and $c_{2}.$ Let $\kappa:=c_{1}+c_{3}.$ It follows from \eqref{eq: ce203} that \begin{equation}\label{eq: ce2014}\|D^{\gamma} (\kappa+\L^{\phi})^{-1}f\|_{\infty}\leq C \|f\|_{2}\quad \forall\; |\gamma|\leq l\end{equation} with $C$ depending only on $c_{1}, c_{2}, c_{3}, n, m.$ By the semigroup property, $$t\|\frac{d}{dt} e^{-t(\L^{\phi} +\kappa)}f\|_{2}\leq C \|f\|_{2}, \;\;\forall \;t>0,$$ where the constant $C$ depends only on $c_{1}, c_{2}, c_{3}.$ Write $$e^{-t(\L^{\phi} +\kappa)}=(\L^{\phi} +\kappa)^{-1}(-1)\frac{d}{dt} e^{-t(\L^{\phi} +\kappa)}.$$ Employing \eqref{eq: ce2014} we derive \begin{equation}\label{eq: ce2015}\|D^{\gamma} e^{-t\L^{\phi} }f\|_{\infty}\leq Ct^{-1}e^{t\kappa} \|f\|_{2}\quad \forall\; t>0.\end{equation} For any fixed $\ez\geq c,$ $L^{\ez\phi} \in \Ez_{2m}(c_{1}, 0, c_{2}, c_{3}\ez^{2m})$ thanks to $\L \in \Ez_{2m}(c_{1}, 0, c_{2}).$ Set $s=\ez^{-1}.$ Then, by scaling, $(L^{\ez\phi})^{s}=s^{2m}T^{-1}L^{\ez\phi}T\in \Ez_{2m}(c_{1}, 0, c_{2}, c_{3})$ and $e^{-L^{\ez\phi} }=Te^{-s^{-2m}(L^{\ez\phi})^{s} }T^{-1}.$  Consequently, applying \eqref{eq: ce2015} with $\L, t$ replaced by $(L^{\ez\phi})^{s}, s^{-2m},$ respectively, we arrive at \begin{equation}\label{eq: ce2030}\begin{aligned}
\|D^{\gamma}e^{-L^{\ez\phi}}f\|_{\infty}&= \|D^{\gamma} (Te^{-s^{-2m}(L^{\ez\phi})^{s}}T^{-1})f\|_{\infty}\\[4pt]
&= s^{-|\gamma|}\|D^{\gamma} e^{-s^{-2m}(L^{\ez\phi})^{s}}T^{-1}f\|_{\infty}\\[4pt]
&\leq Cs^{-|\gamma|} s^{2m} e^{\kappa s^{-2m}}\|T^{-1}f\|_{L^{2}}\\[4pt]
& \leq Cs^{-|\gamma|} s^{2m} e^{\kappa s^{-2m}}\|f\|_{L^{2}} s^{-\frac{n}{2}}\leq C e^{a\ez^{2m}}\ez^{|\gamma|+n/2} \|f\|_{L^{2}}.\end{aligned} 
\end{equation} Apparently, the above arguments apply for $(K_{t}^{\ez\phi})^{*}:=e^{\ez\phi}e^{-t\L^{*}}e^{-\ez\phi}.$ So, \eqref{eq: ce2030} holds for $e^{-(L^{\ez\phi})^{*}}=e^{\ez\phi}e^{-\L^{*}}e^{-\ez\phi}$ when $\gamma=0,$ which implies $$\|e^{-L^{\ez\phi}}f\|_{2}\leq Ce^{a\ez^{2m}}\ez^{n/2}\|f\|_{1}.$$
By this inequality, a duality argument and the semigroup property, we have $$\|D^{\gamma}e^{-L^{\ez\phi}}f\|_{\infty}\leq C e^{2a\ez^{2m}}\ez^{n+|\gamma|}\|f\|_{1},$$ which gives us \eqref{eq: ce2011}. 

\noindent{\textbf{Case 2: 2m=n.}} In the critical case, the Sobolev imbedding of $H^{m}$ into $C^{l, \nu}$ fails. To address this difficulty we use \cite[Lemma 23 Chapter 1]{AT}. Like the argument at the beginning of Case 1, for any $\L \in \Ez_{2m}(c_{1}, \lambda_{0}, c_{2}, c_{3}),$  there exists a $\theta\in (\omega, \frac{\pi}{2})$ such that for any $z\in \Gamma_{\pi-\theta}$ with $|z|=1,$ \begin{equation*}(z+\L+\lambda)^{-1}: H^{-m}\to H^{m}\quad \mbox{ is bounded},\end{equation*} where $\lambda:=c_{1}/2+\lambda_{0}+c_{3}.$ Obviously, $z+\L+\lambda$ ($|z|=1$) can be extended to a bounded operator from $W^{m, p}$ into $W^{-m, p}$ for any $1<p<\infty.$ An application of \cite[Lemma 23 Chapter 1]{AT} implies that $z+\L+\lambda$ is invertible from $W^{m, p}$ onto $W^{-m, p}$ 
for $p$ in a neighborhood of $2$. Choosing $p>2$ in this neighborhood so that $L^{2}\subset W^{-m, p},$ then we can use the Sobolev imbedding ($m-\frac{n}{p}>0$) to deduce that  $$(z+\L+\lambda)^{-1}: L^{2}\to C^{l, \nu}$$ for $l=m-[\frac{n}{p}]-1$ and $0<\nu<1.$ Thus we return to the Case 1. 

All in all, we obtain Gaussian estimates \eqref{eq: ce3002}-\eqref{eq: ce3004} for the heat kernel of $\L$ when $2m\geq n$. Then, we can build the fundamental solution $\Gamma(t, x, s, y)$ of $\H$ by the same way in Theorem \ref{theorem: ce52}. Clearly, $\Gamma(t, x, s, y)$ verifies \eqref{eq: ce53}-\eqref{eq: ce57} with an integer $l$ and some $0<\nu<1$ by Remark \ref{remark: ce101}. The proof is complete.

\hfill$\Box$
\begin{remark}\label{remark: ce4010}  The case $l=0$ in Theorem \ref{theorem: ce1010} is actually stated in \cite[Proposition 28]{AT}. According to the definition of $l$ in the proof of Theorem \ref{theorem: ce1010}, $l+\nu$ may be less than or equal to $m-\frac{n}{2},$ hence Theorem \ref{theorem: ce1010} is not included in \cite{AQ}.
\end{remark}


\section*{Availability of data and material}
 Not applicable.
 
 \section*{Competing interests}
 The author declares that they have no competing interests.

\end{document}